\documentclass[preprint,10pt,authoryear]{elsarticleforarxiv}

\usepackage{amssymb}
\usepackage{amsmath}

\usepackage{mathrsfs}
\usepackage{graphicx}
\usepackage{hyperref}
\usepackage[english]{babel}
\usepackage[T1]{fontenc}
\usepackage[utf8]{inputenc}
\usepackage{enumitem}
\usepackage{mathtools}
\usepackage{amsfonts}
\usepackage[dvipsnames]{xcolor}
\usepackage[nottoc]{tocbibind}
\usepackage{graphicx}
\usepackage{geometry}
\usepackage{dsfont}
\usepackage{fancyhdr}
\usepackage{stix}
\usepackage{multirow}
\usepackage{longtable}
\usepackage{array}
\usepackage{pgfplots}
\usepackage{graphics}
\usepackage{tikz}
\usepackage{pgfkeys}
\usetikzlibrary{arrows}
\usepackage{amsthm}
\usepackage{algorithm}
\usepackage{algpseudocode}
\usepackage{etoolbox}
\usepackage{nicematrix}
\usepackage[many]{tcolorbox}
\usepackage{cleveref}
\usepackage{subcaption}
\usepackage{float}
\usepackage{wrapfig}

\makeatletter
\patchcmd{\upbracefill}{\m@th}{\scriptstyle\m@th}{}{}
\patchcmd{\upbracefill}{$\braceld$}{$\scriptstyle\braceld$}{}{}
\patchcmd{\upbracefill}{\bracelu}{\bracelu\mkern-1mu}{}{}
\patchcmd{\upbracefill}{\hfill\braceru}{\hfill\mkern-1mu\braceru}{}{}
\makeatother

\newtheorem{theorem}{Theorem}
\newtheorem{lemma}{Lemma}
\newtheorem{conjecture}{Conjecture}
\newtheorem{observation}{Observation}
\newtheorem{proposition}{Proposition}
\newtheorem{definition}{Definition}
\newtheorem{remark}{Remark}
\newtheorem{remarks}[remark]{Remark}
\newtheorem*{property}{Property}

\newtcolorbox{Corner}{text width=(\linewidth-\parindent),colback=white,colframe=white,enhanced,overlay={
\draw[line width=1pt] ([xshift=2mm,yshift=-1mm]frame.north west) -- ([yshift=-1mm]frame.north west) -- ([yshift=1mm]frame.south west) -- ([xshift=2mm,yshift=1mm]frame.south west) ;}} 

\newenvironment{Theorem}{\begin{Corner}
\begin{theorem}}{\end{theorem}
\end{Corner}}
\newenvironment{Lemma}{\begin{Corner}
\begin{lemma}}{\end{lemma}
\end{Corner}}
\newenvironment{Conjecture}{\begin{Corner}
\begin{conjecture}}{\end{conjecture}
\end{Corner}}

\newenvironment{Proposition}{\begin{Corner}
\begin{proposition}}{\end{proposition}
\end{Corner}}
\newenvironment{Definition}{\begin{Corner}
\begin{definition}}{\end{definition}
\end{Corner}}
\newenvironment{Remark}{\begin{Corner}
\begin{remark}}{\end{remark}
\end{Corner}}
\newenvironment{Remarks}[1][]{\begin{Corner}
\begin{remarks}[#1]\begin{itemize}[label=$\bullet$,leftmargin=0mm]\item[]}{\end{itemize}\end{remarks}
\end{Corner}}

\newenvironment{Property}{\begin{Corner}\begin{property}}{\end{property}\end{Corner}}

\definecolor{CouleurLP}{RGB}{226,111,70}
\definecolor{CouleurSDPNaif}{RGB}{61,164,77}
\definecolor{CouleurSDPQuadraticInequalities}{RGB}{0,80,230}
\definecolor{CouleurSDPMISC}{RGB}{150,32,72}

\newcommand{\exposant}[1]{\ensuremath{\cdot\text{\footnotesize 10}^{\hspace*{1pt}\text{#1}}}}
\newcommand{\exposantvanilla}[1]{\ensuremath{\text{\footnotesize 10}^{\hspace*{1pt}\text{#1}}}}
\newcommand\Tstrut{\rule{0pt}{2.6ex}}
\newcommand\Bstrut{\rule[-0.9ex]{0pt}{0pt}}

\DeclareMathOperator{\NP}{\mathscr{N}\hspace*{-3pt}\mathscr{P}}
\DeclareMathOperator{\tr}{tr}
\DeclareMathOperator{\N}{\mathbb{N}}

\DeclareMathOperator{\R}{\mathbb{R}}
\DeclareMathOperator{\SDP}{\mathbb{S}}
\DeclareMathOperator{\diag}{diag}
\DeclareMathOperator{\Diag}{Diag}
\DeclareMathOperator{\Opt}{Opt}
\DeclareMathOperator{\conv}{conv}
\DeclareMathOperator{\argmin}{argmin}

\newcounter{FigureCounter}
\newcounter{EquationCounter}
\addto\captionsenglish{}

\pgfplotsset{compat=1.18}
\usepackage[fixed]{fontawesome5}

\definecolor{CouleurCite}{RGB}{37,52,190}

\newtagform{noparen}{\color{CouleurCite}}{}
\usetagform{noparen}
\newcommand{\refequation}[1]{\textcolor{CouleurCite}{\ref{#1}}}
\newcommand{\citeequation}[1]{\textcolor{CouleurCite}{\citet{#1}}}

\crefname{figure}{\textsc{Figure}}{\textsc{Figures}}

\begin{document}

\begin{frontmatter}



\title{Knapsack with compactness: a semidefinite approach}


\author{Hubert Villuendas\corref{cor1}} 
\cortext[cor1]{Corresponding author\\
\emph{Email adresses}: \texttt{$\left\lbrace\,\text{hubert.villuendas}\,,\,\text{mathieu.besancon}\,,\,\text{jerome.malick}\,\right\rbrace$@univ-grenoble-alpes.fr}}

\affiliation{organization={University Grenoble Alpes, CNRS, Inria, Grenoble INP, LJK, LIG},
            country={France}}
            
\author{Mathieu Besançon} 

\affiliation{organization={University Grenoble Alpes, Inria, CNRS, LIG},
            country={France}}

\author{Jérôme Malick} 

\affiliation{organization={University Grenoble Alpes, CNRS, Grenoble INP, LJK},
            country={France}}
            
\begin{abstract}
The min-knapsack problem with compactness constraints extends the classical knapsack problem, in the case of ordered items, by introducing a restriction ensuring that they cannot be too far apart. This problem has applications in statistics, particularly in the detection of change-points in time series. In this paper, we propose a semidefinite programming approach for this problem, incorporating compactness in constraints or in objective. We study and compare the different relaxations, and argue that our method provides high-quality heuristics and tight bounds. In particular, the single hyperparameter of our penalized semidefinite models naturally balances the trade-off between compactness and accuracy of the computed solutions. Numerical experiments illustrate, on the hardest instances, the effectiveness and versatility of our approach compared to the existing mixed-integer programming formulation.
\end{abstract}




\end{frontmatter}



\section{Introduction}

\subsection{Knapsack and compactness}

We consider the min-knapsack problem with a "compactness" constraint. Given $n$ items, each with a weight $w_i\geq 0$ and a cost $c_i\geq 0$, the min-knapsack consists in selecting  a subset of elements of minimum cost and whose weight is at least a certain threshold. This problem formulates as
\begin{equation}
\left[\begin{array}{rl}
\text{minimize} & c^\top x\\
\text{ subject to} & w^\top x\geq q\\
& x\in\left\lbrace 0,1\right\rbrace^n
\end{array}\right.\tag{(mKP)}\label{eq: mKP sans compactness}
\end{equation}
which is equivalent to the classical maximization knapsack \citeequation{csirik1991heuristics}. Indeed, \refequation{eq: mKP sans compactness} is related to its "complement" instance:
\begin{equation}
\left[\begin{array}{rl}
\text{minimize} & c^\top x\\
\text{subject to} & w^\top x\geq q\\
& x\in\left\lbrace 0,1\right\rbrace^n
\end{array}\right.\qquad\Leftrightarrow\qquad \left[\begin{array}{rl}
\text{maximize} & c^\top z\\
\text{subject to} & w^\top z \leq\displaystyle\sum\limits^n_{i=1}w_i-q\\
& z\in\left\lbrace 0,1\right\rbrace^n.
\end{array}\right.\tag{$\left(\theEquationCounter\right)$}\label{eq: Explicitation de la bijection entre min-knapsack et knapscak classique}
\end{equation}\addtocounter{EquationCounter}{1}\null
and we have an explicit bijection $x\mapsto z:=\mathbf{1}-x$ between the solutions of the two problems.

Furthermore, we consider a unidimensional ordering of the items and enforce compactness: we ensure that selected elements remain close to one another within this sequence, thereby preserving a degree of structural coherence in the solution. The compactness of the solution is achieved in \citeequation{santini2024min} by enforcing a \emph{compactness constraint} added to \refequation{eq: mKP sans compactness}: for a given integer $\Delta\in\N$, the distance $\vert i-j\vert$ separating two consecutive selected objects $i,j$ must not exceed $\Delta$. For example, with $\Delta=2$, this constraint is not verified by items 1 and 4 on \cref{Fig: non-compact selection}; adding item 2 yields the compact selection proposed on \cref{Fig: compact selection}.
\begin{figure}[H]
   \begin{minipage}{0.45\linewidth}
       \begin{center}
       \begin{tikzpicture}[line cap=round,line join=round,>=triangle 45,x=1cm,y=1cm]
       \clip(-3.2,-1) rectangle (3.2,1);
       \draw (-3,0.1) node[anchor=south] {$1$};
       \draw (0,0.1) node[anchor=south] {$4$};
       \draw (1,0.1) node[anchor=south] {$5$};
       \draw (3,0.1) node[anchor=south] {$7$};
       \draw [<->,> = stealth,line width=1pt] (-2.8,0) -- (-0.2,0);
       \draw [<->,> = stealth,line width=1pt] (0.2,0) -- (0.8,0);
       \draw [<->,> = stealth,line width=1pt] (1.2,0) -- (2.8,0);
       \begin{scriptsize}
       \draw [fill=black] (-3,0) circle (2.5pt);
       \draw [fill=white] (-2,0) circle (2pt);
       \draw [fill=white] (-1,0) circle (2pt);
       \draw [fill=black] (0,0) circle (2.5pt);
       \draw [fill=black] (1,0) circle (2.5pt);
       \draw [fill=white] (2,0) circle (2pt);
       \draw [fill=black] (3,0) circle (2.5pt);
       \draw (-1.5,-0.1) node[anchor=north] {\textcolor{red}{$3>\Delta$}};
       \draw (0.5,-0.1) node[anchor=north] {\textcolor{green!70!black}{$1\leq\Delta$}};
       \draw (2,-0.1) node[anchor=north] {\textcolor{green!70!black}{$2\leq\Delta$}};
       \end{scriptsize}
       \end{tikzpicture}
       \vspace*{-8mm}\caption{a \textcolor{red}{\textbf{non}-compact} selection with $\Delta =2$ and $n=7$.}
       \label{Fig: non-compact selection}       
       \end{center}
   \end{minipage}\hfill
   \begin{minipage}{0.45\linewidth}
       \begin{center}
           \begin{tikzpicture}[line cap=round,line join=round,>=triangle 45,x=1cm,y=1cm]
           \clip(-3.2,-1) rectangle (3.2,1);
           \draw (-3,0.1) node[anchor=south] {$1$};
           \draw (-2,0.1) node[anchor=south] {$2$};
           \draw (0,0.1) node[anchor=south] {$4$};
           \draw (1,0.1) node[anchor=south] {$5$};
           \draw (3,0.1) node[anchor=south] {$7$};
           \draw [<->,> = stealth,line width=1pt] (-2.8,0) -- (-2.2,0);
           \draw [<->,> = stealth,line width=1pt] (-1.8,0) -- (-0.2,0);
           \draw [<->,> = stealth,line width=1pt] (0.2,0) -- (0.8,0);
           \draw [<->,> = stealth,line width=1pt] (1.2,0) -- (2.8,0);
           \begin{scriptsize}
           \draw [fill=black] (-3,0) circle (2.5pt);
           \draw [fill=black] (-2,0) circle (2.5pt);
           \draw [fill=white] (-1,0) circle (2pt);
           \draw [fill=black] (0,0) circle (2.5pt);
           \draw [fill=black] (1,0) circle (2.5pt);
           \draw [fill=white] (2,0) circle (2pt);
           \draw [fill=black] (3,0) circle (2.5pt);
           \draw (-2.5,-0.1) node[anchor=north] {\textcolor{green!70!black}{$1\leq\Delta$}};
           \draw (-1,-0.1) node[anchor=north] {\textcolor{green!70!black}{$2\leq\Delta$}};
           \draw (0.5,-0.1) node[anchor=north] {\textcolor{green!70!black}{$1\leq\Delta$}};
           \draw (2,-0.1) node[anchor=north] {\textcolor{green!70!black}{$2\leq\Delta$}};
           \end{scriptsize}
           \end{tikzpicture}
       \vspace*{-8mm}\caption{a \textcolor{green!70!black}{compact} selection with $\Delta =2$ and $n=7$.}
       \label{Fig: compact selection} 
       \end{center}
       \end{minipage}
\end{figure}

This constraint can be modeled by requiring that, for any pair of indices $i$ and $j$ distant enough, if both items $i$ and $j$ are selected, then at least one item between them must also be included in the selection. Assuming without loss of generality that $j > i$, this condition can be expressed in two ways: with a linear or a quadratic inequality
\begin{center}
\begin{minipage}{0.48\linewidth}
\begin{equation}
x_i+x_j-1\leq\displaystyle\sum\limits_{k=i+1}^{j-1}x_k\tag{$\left(\theEquationCounter\text{a}\right)$}\label{eq: inegalite de compacite version lineaire}
\end{equation}
\end{minipage}\hfill\begin{minipage}{0.48\linewidth}
\begin{equation}
x_ix_j\leq\displaystyle\sum\limits_{k=i+1}^{j-1}x_k.\tag{$\left(\theEquationCounter\text{b}\right)$}\label{eq: inegalite de compacite version quadratique}
\end{equation}\addtocounter{EquationCounter}{1}\null
\end{minipage}
\end{center}

Existing work \citeequation{santini2024min,curebal2024fixed} about compactness enfore it by adding the constraint \refequation{eq: inegalite de compacite version lineaire}, or a strengthened form of \refequation{eq: inegalite de compacite version lineaire}, to the mixed-integer formulation \refequation{eq: original formulation of the mKPC}. Thus, \citeequation{santini2024min} studies the following linear binary formulation of the min-knapsack problem with compactness constraints:
\begin{equation}
\left[\begin{array}{rlc}
\text{minimize} & c^\top x &\\
\text{subject to} & w^\top x\geq q &\\
& \forall i,j\in[\![n]\!]\text{ with }j-i>\Delta,\quad \left\lfloor\dfrac{j-i-1}{\Delta}\right\rfloor\left(x_i+x_j-1\right)\leq\displaystyle\sum\limits_{k=i+1}^{j-1}x_k &\quad (\bigstar)\\
& x\in\left\lbrace 0,1\right\rbrace^n &
\end{array}\right.\tag{(mKPC)}\label{eq: original formulation of the mKPC}
\end{equation}
The constraints  \hyperref[eq: original formulation of the mKPC]{$(\bigstar)$} ensure that for all $i,j\in [\![n]\!]$, if items $i$ and $j$ are selected and are sufficiently far apart, then at least $\left\lfloor (j-i-1)/\Delta\right\rfloor$ items between $i$ and $j$ must also be selected. We denote by $\textcolor{CouleurCite}{\left(\text{mKPC}\right)_{\text{LP}}}\label{eq: original formulation of the mKPC - LP version}$ its linear relaxation, replacing the constraint $x\in\left\lbrace 0,1\right\rbrace^n$ by $x\in[0,1]^n$.

As a variant of the knapsack problem, \refequation{eq: original formulation of the mKPC} is inherently $\NP$-hard; it is moreover significally more challenging to solve than the classical version. Indeed, while the classical knapsack problem is often relatively easy to solve in practice, to the point that it is often considered one of the easiest $\NP$-hard problems \citeequation{pisinger2005hard}, the addition of compactness constraints \hyperref[eq: original formulation of the mKPC]{$(\bigstar)$} significantly increases computational complexity. Moreover, problem difficulty varies considerably depending on instance characteristics, with certain instances proving substantially harder to solve using off-the-shelf solvers, due to the fact that the linear relaxation \hyperref[eq: original formulation of the mKPC - LP version]{$\textcolor{CouleurCite}{\left(\text{mKPC}\right)_{\text{LP}}}$} often yields a large number of fractional variables with values close to $1/2$, resulting in significantly more branching during the exploration of the branch~{\&}~bound tree, as observed by \citeequation{santini2024min}. Then \citeequation{curebal2024fixed} proposes to combine a heuristic with the integer programming approach.

In this work, we investigate an alternative approach: exploiting the quadratic formulation of the compactness inequality \refequation{eq: inegalite de compacite version quadratique}, we consider the min-knapsack with compactness as a binary problem with quadratic constraints and propose strong relaxations based on semidefinite programming, aiming both to derive tight bounds and to guide the development of effective heuristics. In particular, we seek to avoid the issue of poor fractionality in the computed relaxed solution, which poses a significant challenge for exact methods. We provide more details on our approach, our results, and the outline of the paper, in \cref{subsec: contributions and outline}. Before this, in the next section, we give an example of how this compactness appears in an application.


\subsection{Example of application}\label{sec: motivation}

The min-knapsack with compactness captures applications in statistics, particularly in the detection of regime changes in time series. In this section, we briefly present this application, refering to \citeequation{santini2024min,pisinger2005hard,aminikhanghahi2017survey} for more details.

Given a sequence of measurements collected at regular intervals, describing the behaviour of a system, the aim of change point detection is to identify abrupt breaks in a time series, signaling significant changes in the dynamics or distribution of the data as a result of external events. Change detection in time series has numerous practical applications \citeequation{aminikhanghahi2017survey}. In climatology, it can be used to identify the effects of climate change. In medicine, continuous monitoring of a patient's state of health relies on real-time analysis of physiological variables (such as heart rate, electroencephalogram), see \emph{e.g.} \citeequation{aminikhanghahi2017survey} the recent survey on time series and their applications. 

\begin{figure}[H]
   \begin{minipage}{0.45\linewidth}
       \centering
       \includegraphics[width=\linewidth,keepaspectratio]{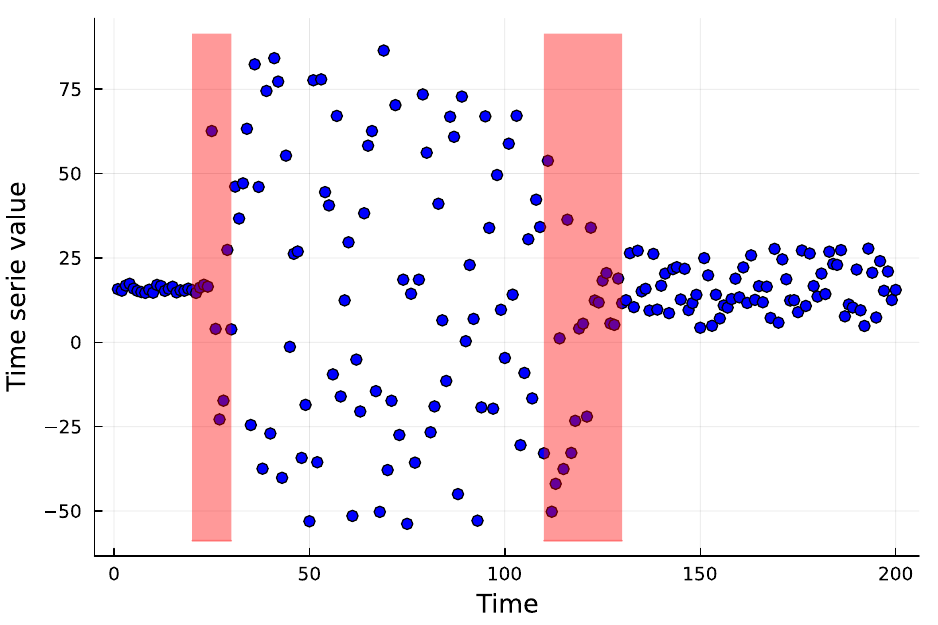}
       \vspace*{-8mm}\caption{a time serie and its possible changes in variance.}
       \label{Fig: a time serie and its possible changes in variance}
   \end{minipage}\hfill
   \begin{minipage}{0.45\linewidth}
       \centering
       \includegraphics[width=\linewidth,keepaspectratio]{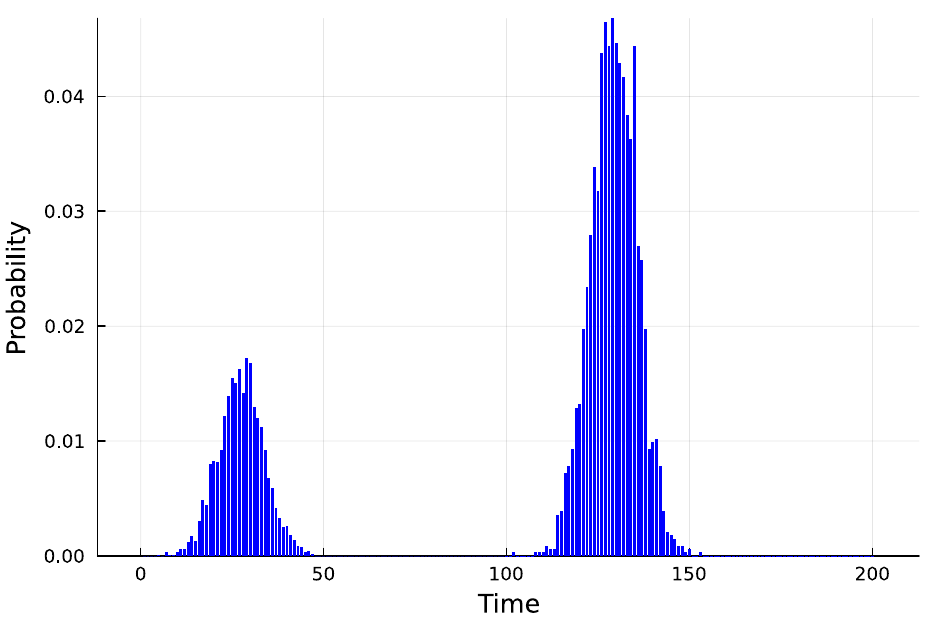}
       \vspace*{-8mm}\caption{probabilities associated to each point of being the first change point of the time serie.}
       \label{Fig: probabilites associated to each point}
   \end{minipage}
\end{figure}

For detecting variance variations in Gaussian time series, \citeequation{cappello2022bayesian} proposed an iterative procedure that identifies a change point at each step. However, assigning a single point in time to each change can be difficult to interpret because of the uncertainty on precise moment at which the variation occurs; see \emph{e.g.} \cref{Fig: a time serie and its possible changes in variance}.

\begin{figure}[H]
   \begin{minipage}{0.45\linewidth}
       \centering
       \includegraphics[width=\linewidth,keepaspectratio]{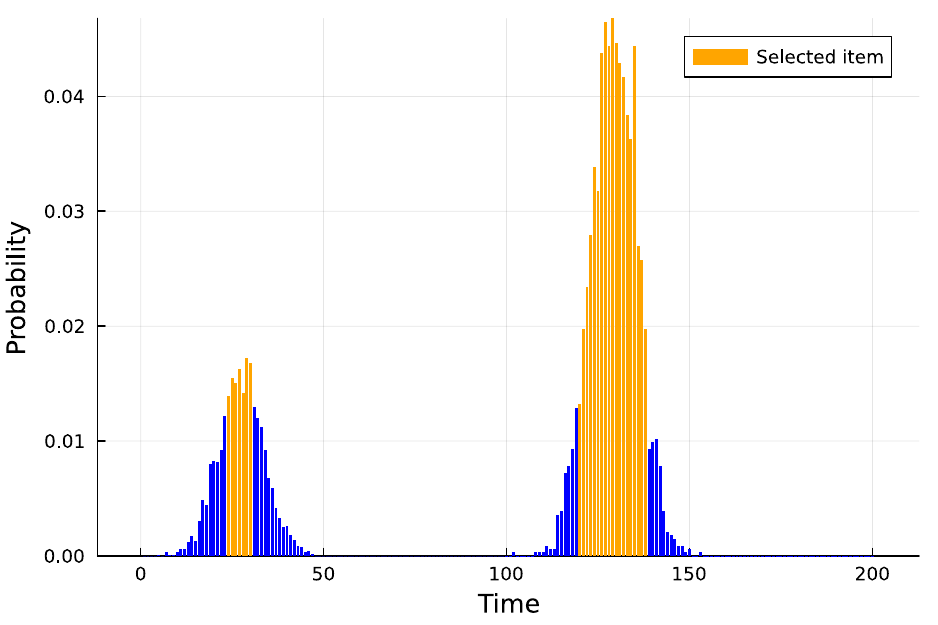}
       \vspace*{-8mm}\caption{selected items for the \refequation{eq: mKP sans compactness}, building a \textcolor{red}{\textbf{non}-compact} $0.75$-credible set, with a compactness parameter $\Delta = 1$.}
       \label{Fig: selected items with the (mKP)}
   \end{minipage}\hfill
   \begin{minipage}{0.45\linewidth}
       \centering
       \includegraphics[width=\linewidth,keepaspectratio]{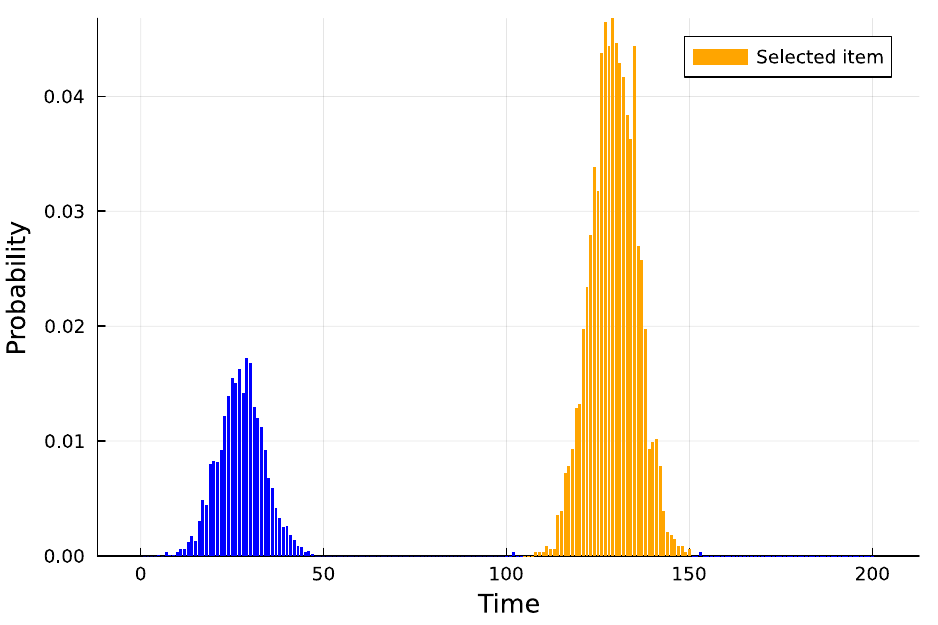}
       \vspace*{-8mm}\caption{selected items for the \refequation{eq: original formulation of the mKPC}, and the corresponding \textcolor{green!70!black}{compact} $0.75$-credible set, with a compactness parameter $\Delta = 1$.}
       \label{Fig: selected items with the (mKPC)}
   \end{minipage}
\end{figure}

To overcome this problem, the method builds a discrete probability distribution over the time points, associating with each one a probability that it constitutes a point of change. \cref{Fig: probabilites associated to each point} displays an example: the bars on the graph represent the probabilities assigned to each instant in the time series. The method then identifies a subset whose associated probabilities sum to at least $q$, while minimizing the number of selected elements. To account for varying measurement accuracy, each point $i$ is assigned a cost $c_i \geq 0$ reflecting its inaccuracy, and we aim to minimize the total imprecision. This leads to a standard min-knapsack formulation. However, solving the problem as \refequation{eq: mKP sans compactness} may result in sets that mix multiple change points, as illustrated in \cref{Fig: selected items with the (mKP)}. Since each selected set is meant to capture a single change point, the selection should ideally be compact. This motivates the compactness constraint in \refequation{eq: original formulation of the mKPC}, whose effect is visualized in \cref{Fig: selected items with the (mKPC)}.

\subsection{Contributions and outline}\label{subsec: contributions and outline}

In this paper, we introduce a semidefinite reformulation of the min-knapsack problem with compactness, leading to two distinct semidefinite relaxations: the first relaxation explicitly incorporates the compactness constraint within the formulation, building on the approach on \citeequation{santini2024min}, while the second indirectly promotes compactness through a penalized objective function. To improve the quality of these relaxations, we investigate various strengthening techniques and propose two families of valid inequalities: one leveraging semidefinite programming principles and the other exploiting the combinatorial structure of the problem. In particular, we adapt key concepts from classical knapsack literature to our framework and introduce an algorithm capable of generating effective cutting planes in pseudo-polynomial time. We then analyze the computational efficiency and solution quality of these approaches, with a particular focus on the challenging instances of \citeequation{santini2024min}. Additionally, we demonstrate how these elements can be combined to design an efficient heuristic that produces high-quality solutions in terms of fractionality, while allowing for a tunable trade-off between accuracy of the solution set and compactness.

This paper is organized as follows: in \cref{sec: Modeling the min-knapsack with compactness as a Semidefinite Program}, we present our semidefinite formulations of the min-knapsack problem with compactness constraints and its penalized variant. Section~\ref{sec: Strenghtening the semidefinite relaxation} details strengthening techniques, involving valid inequalities from the theory of semidefinite programming and a separation procedure that rely on the combinatorial properties of the problem. Finally, \cref{sec: computationnal results} reports the results of our numerical experiments, and we conclude with discussions on future research directions. Two appendices complement the paper, for additional information, including a counter-example on a curious observation in \cref{sec: Modeling the min-knapsack with compactness as a Semidefinite Program}.

\section{Modeling as a semidefinite program}\label{sec: Modeling the min-knapsack with compactness as a Semidefinite Program}

In this section, we develop our semidefinite programming approach for the min-knapsack problem with compactness constraints. We first write the exact SDP reformulation by lifting the binary variables to a rank one matrix variable, thereby establishing a formal equivalence with the original formulation. We then present a first semidefinite relaxation obtained from this reformulation, that we compare to the linear relaxation \hyperref[eq: original formulation of the mKPC - LP version]{$\textcolor{CouleurCite}{\left(\text{mKPC}\right)_{\text{LP}}}$}. Finally, we propose a penalized SDP formulation that integrates the compactness requirement directly into the objective function, offering a flexible trade-off between minimizing over the sum of the costs and ensuring a structural compactness.

\subsection{Equivalent SDP reformulations}

In the context of a $\left\lbrace 0,1\right\rbrace$-formulation of the min-knapsack problem with compactness constraint, we apply the usual "recipe" \citeequation{poljak1995recipe} and get a semidefinite reformulation by substituting a matrix variable $\mathbf{X}=xx^\top$ into the original problem. For the sake of completness, we formalize the result and include a direct proof of the equivalence between the two formulation. We will note $\SDP^n$ the set of symmetric matrix of size $n\times n$ and $\SDP^n_+$ the set of symmetric positive semidefinite matrices of size $n\times n$.

\begin{Proposition}\label{prop: equivalence entre (mKPC) et la formulation SDP avec rang} Let $\mathcal{X}\subseteq\SDP^n$ denote the set of all symmetric matrices $\mathbf{X}$ that verify the linear inequalities
\begin{equation}
\mathcal{X}=\left\lbrace\mathbf{X}\in\SDP^n\,\left\vert\,w^\top\diag(\mathbf{X})\geq q\text{ and }\forall i,j\in[\![n]\!],j-i>\Delta, \left\lfloor\dfrac{j-i-1}{\Delta}\right\rfloor\mathbf{X}_{ij}\leq\sum^{j-1}_{k=i+1}\mathbf{X}_{kk}\right.\right\rbrace\tag{$\left(\theEquationCounter\right)$}\label{eq: definition du polytope de la proposition 1}
\end{equation}\addtocounter{EquationCounter}{1}\null
where $\diag:\SDP^n\rightarrow\R^n$ is the operator taking a matrix $\mathbf{X}$ and associates the vector of its diagonal entries.
Then both \emph{\refequation{eq: original formulation of the mKPC}} and
\begin{equation}
\min_{\mathbf{X}\in\SDP^n_+}\left\lbrace c^\top\diag\left(\mathbf{X}\right)\,\left\vert\,\mathbf{X}\in\mathcal{X}\cap\left\lbrace 0,1\right\rbrace^{n\times n},\,\mathrm{rank}\left(\mathbf{X}\right)=1\right.\right\rbrace\tag{$\left(\theEquationCounter\right)$}\label{eq: probleme SDP equivalent au probleme en 0,1}
\end{equation}\addtocounter{EquationCounter}{1}\null
yield the same optimal solution. More precisely, if $x^\ast$ is an optimal solution of \emph{\refequation{eq: original formulation of the mKPC}}, then $\mathbf{X}^\ast= x^\ast{x^\ast}^\top$ is an optimal solution of \refequation{eq: probleme SDP equivalent au probleme en 0,1}; conversely, if $\mathbf{X}^\ast$ is an optimal solution of \refequation{eq: probleme SDP equivalent au probleme en 0,1}, then $x^\ast=\diag\left(\mathbf{X}^\ast\right)$ is an optimal solution of \emph{\refequation{eq: original formulation of the mKPC}}.
\end{Proposition}
\begin{proof}
Let $x^\ast$ be an optimal solution of \refequation{eq: original formulation of the mKPC}. Then we set $\mathbf{X}=x^\ast{x^\ast}^\top$, and since $x^\ast$ is a $\left\lbrace 0,1\right\rbrace$-vector, we have $\mathbf{X}_{kk}=\left(x^\ast_k\right)^2=x^\ast_k$ for all $k\in[\![n]\!]$. As a consequence:
$$\sum^n_{i=1}w_i\mathbf{X}_{ii}=\sum^n_{i=1}w_i\left(x^\ast_i\right)^2=\sum^n_{i=1}w_ix^\ast_i\geq q.$$
Moreover, if $i,j\in[\![n]\!]$ are such that $j-i>\Delta$, two cases occur: either at least one of the variables is zero, and then $\mathbf{X}_{ij}=x^\ast_ix^\ast_j=0$ so that the inequality
$$\left\lfloor\dfrac{j-i-1}{\Delta}\right\rfloor\mathbf{X}_{ij}=0\leq\sum^{j-1}_{k=i+1}\mathbf{X}_{kk}$$
is indeed verified; either $x^\ast_i=x^\ast_j=1$ and we then have $\mathbf{X}_{ij}=1$ and $x^\ast_i+x^\ast_j-1=1$, leading to
$$\left\lfloor\dfrac{j-i-1}{\Delta}\right\rfloor\mathbf{X}_{ij}=\left\lfloor\dfrac{j-i-1}{\Delta}\right\rfloor\left(x^\ast_i+x^\ast_j-1\right)\leq\sum^{j-1}_{k=i+1}x^\ast_k=\sum^{j-1}_{k=i+1}\mathbf{X}_{kk}$$
where the inequality comes from the fact that $x^\ast$ is a solution of \refequation{eq: original formulation of the mKPC}. Thus $\mathbf{X}\in\mathcal{X}$, has rank one and its coefficients in $\left\lbrace 0,1\right\rbrace$, so $\mathbf{X}$ is a solution of \refequation{eq: probleme SDP equivalent au probleme en 0,1}, with the same objective value $c^\top\diag\left(\mathbf{X}\right)=c^\top x^\ast$. In particular, the optimal value of \refequation{eq: probleme SDP equivalent au probleme en 0,1} is less than the optimal value of the \refequation{eq: original formulation of the mKPC}.

Conversely, if $\mathbf{X}^\ast$ is an optimal solution of \refequation{eq: probleme SDP equivalent au probleme en 0,1}, then $\mathbf{X}^\ast$ is a positive semidefinite matrix of rank one with coefficients in $\left\lbrace 0,1\right\rbrace$, so there is a vector $x\in\left\lbrace 0,1\right\rbrace^n$ such that $\mathbf{X}^\ast=xx^\top$. Applying the fact that for all $k\in[\![n]\!]$, $\mathbf{X}^\ast_{kk}=x_k^2=x_k$, we have
$$w^\top x=\sum^n_{i=1}w_ix_k=\sum^n_{i=1}w_i\mathbf{X}^\ast_{kk}\geq q$$
and for all $i,j\in[\![n]\!]$, we have $\left(1-x_i\right)\left(1-x_j\right)\geq 0$, which gives us the inequality $x_ix_j\geq x_i+x_j-1$, so if $j-i>\Delta$:
$$\left\lfloor\dfrac{j-i-1}{\Delta}\right\rfloor\left(x_i+x_j-1\right)\leq\left\lfloor\dfrac{j-i-1}{\Delta}\right\rfloor\underbrace{x_ix_j}_{=\mathbf{X}^\ast_{ij}}\leq\sum^{j-1}_{k=i+1}x_k.$$
So $x$ is indeed a solution of the \refequation{eq: original formulation of the mKPC} with value $c^\top x=w_1\mathbf{X}^\ast_{11}+\dots+w_n\mathbf{X}^\ast_{nn}$, giving us the equality between the value of \refequation{eq: probleme SDP equivalent au probleme en 0,1} and the optimal value of \refequation{eq: original formulation of the mKPC}, and we have the equality $\mathbf{X}^\ast=x^\ast{x^\ast}^\top$ and $\diag\left(\mathbf{X}^\ast\right)=x^\ast$.
\end{proof}

This equivalent reformulation of \refequation{eq: original formulation of the mKPC} shows non-convexity in both the rank constraint and the $\left\lbrace 0,1\right\rbrace$ coefficient constraints. We use a classical result about $\left\lbrace 0,1\right\rbrace$ positive semidefinite matrices, recalled in \ref{app: schur complement lemma}.
 
\begin{Proposition}\label{prop: formulation equivalente avec la contrainte conique} Let $\mathcal{X}\subseteq\SDP^n$ be the polytope defined at \refequation{eq: definition du polytope de la proposition 1} in \textbf{\emph{Proposition \ref{prop: equivalence entre (mKPC) et la formulation SDP avec rang}}}. Then both \emph{\refequation{eq: original formulation of the mKPC}} and
\begin{equation}
\min_{\mathbf{X}\in\SDP^n_+}\left\lbrace c^\top\diag\left(\mathbf{X}\right)\,\left\vert\,\mathbf{X}\in\mathcal{X},\,\begin{pmatrix}
   1 & \diag(\mathbf{X})^\top\\
   \diag(\mathbf{X}) & \mathbf{X}
\end{pmatrix}\succeq 0\,,\,\mathrm{rank}\left(\mathbf{X}\right)=1\right.\right\rbrace\tag{$\left(\theEquationCounter\right)$}\label{eq: equivalence entre (mKPC) et la formulation SDP avec la matrice par blocs}
\end{equation}\addtocounter{EquationCounter}{1}\null
yield equivalent optimal solution: $x^\ast=\diag\left(\mathbf{X}^\ast\right)$ and $\mathbf{X}^\ast={x^\ast}{x^\ast}^\top$.
\end{Proposition}
\begin{proof}
By \textbf{Proposition \ref{prop: equivalence entre (mKPC) et la formulation SDP avec rang}}, if $x^\ast$ is an optimal solution of \refequation{eq: original formulation of the mKPC}, then $\mathbf{X}={x^\ast}{x^\ast}^\top$ is an optimal solution of \refequation{eq: probleme SDP equivalent au probleme en 0,1}, and since $x^\ast$ is a $\left\lbrace 0,1\right\rbrace$-vector, we have
\begin{equation*}
 \begin{pmatrix}
       1 & \diag\left(\mathbf{X}\right)^\top\\
       \diag\left(\mathbf{X}\right)^\top & \mathbf{X}
\end{pmatrix} = \begin{pNiceArray}{cc}[margin]
1 & {x_1^\ast}^2\,\,\cdots\,\,{x_n^\ast}^2\\
{x_1^\ast}^2 & \\
\vdots &\Block{1-1}<\Large>{\mathbf{X}} \\
{x_n^\ast}^2 & \\
\end{pNiceArray} = \begin{pNiceArray}{cc}[margin]
1 & {x_1^\ast}\,\,\cdots\,\,{x_n^\ast}\\
{x_1^\ast} & \\
\vdots &\Block{1-1}<\Large>{{x^\ast}{x^\ast}^\top} \\
{x_n^\ast} & \\
\end{pNiceArray} = \begin{pmatrix}
   1\\
   {x^\ast}
\end{pmatrix}\begin{pmatrix}
   1\\
   {x^\ast}
\end{pmatrix}^\top\succeq 0 
\end{equation*}
hence $\mathbf{X}$ is indeed a solution of \refequation{eq: equivalence entre (mKPC) et la formulation SDP avec la matrice par blocs}. In particular, \textbf{Proposition \ref{prop: equivalence entre (mKPC) et la formulation SDP avec rang}} states that $\mathbf{X}$ is an optimal solution of \refequation{eq: probleme SDP equivalent au probleme en 0,1}. Since \refequation{eq: probleme SDP equivalent au probleme en 0,1} and \refequation{eq: equivalence entre (mKPC) et la formulation SDP avec la matrice par blocs} both have the same objective function, we have $\Opt\text{\refequation{eq: equivalence entre (mKPC) et la formulation SDP avec la matrice par blocs}}\leq\Opt\text{\refequation{eq: probleme SDP equivalent au probleme en 0,1}}$.

Conversely, if $\mathbf{X}^\ast$ is an optimal solution of \refequation{eq: equivalence entre (mKPC) et la formulation SDP avec la matrice par blocs}, then \textbf{Theorem \ref{thm: avec la matrice par bloc, equivalence entre X est de rang 1 et X est en 0,1}} (see \ref{app: schur complement lemma}) ensures that $\mathbf{X}^\ast$ is a $\left\lbrace 0,1\right\rbrace$ matrix with $\mathrm{rank}\left(\mathbf{X}^\ast\right)=1$, thus $\mathbf{X}^\ast$ is a solution of \refequation{eq: probleme SDP equivalent au probleme en 0,1}. Hence $\Opt\text{\refequation{eq: probleme SDP equivalent au probleme en 0,1}}\leq\Opt\text{\refequation{eq: equivalence entre (mKPC) et la formulation SDP avec la matrice par blocs}}$, which implies that problems \refequation{eq: probleme SDP equivalent au probleme en 0,1} and \refequation{eq: equivalence entre (mKPC) et la formulation SDP avec la matrice par blocs} have the same objective value, and $\mathbf{X}^\ast$ is an optimal solution of \refequation{eq: probleme SDP equivalent au probleme en 0,1}. By \textbf{Proposition \ref{prop: equivalence entre (mKPC) et la formulation SDP avec rang}}, $x=\diag\left(\mathbf{X}^\ast\right)$ is an optimal solution of the \refequation{eq: original formulation of the mKPC}.
\end{proof}

\subsection{A first SDP relaxation}\label{subsec: a first sdp relaxation}

\begin{figure}[H]
\centering
   \begin{minipage}{0.45\linewidth}
       \centering
       \includegraphics[width=\linewidth,keepaspectratio]{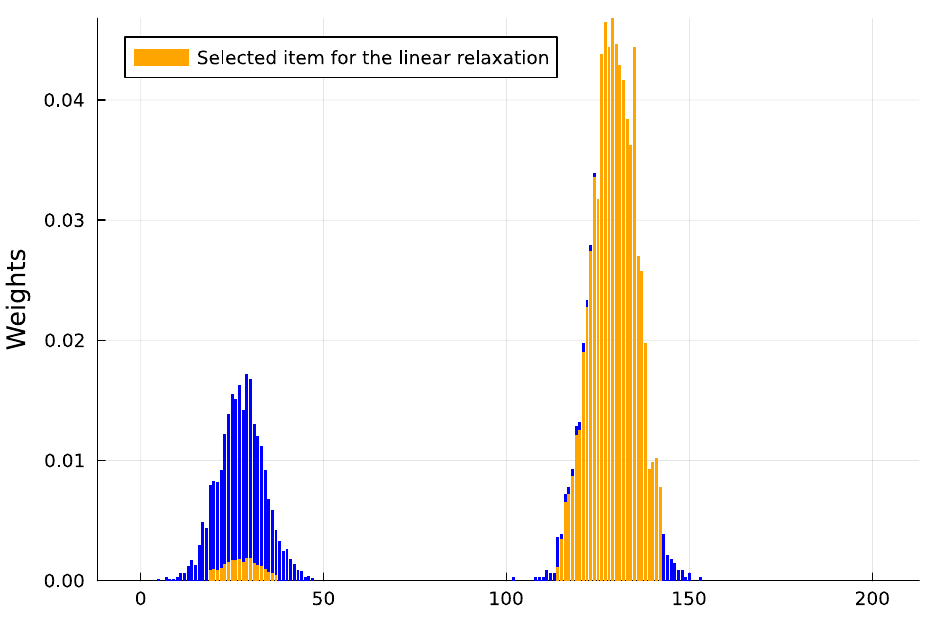}
       \vspace*{-8mm}\caption{selected items for the linear relaxation of \refequation{eq: original formulation of the mKPC}.}
       \label{Fig: selected items with the LP (mKPC)}
   \end{minipage}\hfill
   \begin{minipage}{0.45\linewidth}
       \centering
       \includegraphics[width=\linewidth,keepaspectratio]{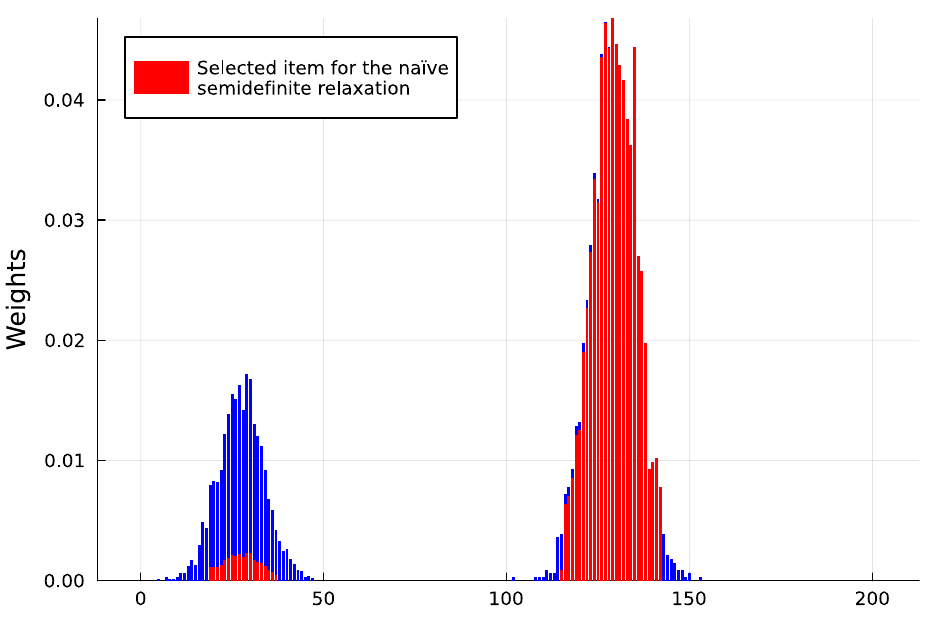}
       \vspace*{-8mm}\caption{selected items for the naïve semidefinite relaxation \refequation{eq: naive semidefinite relaxation}.}
       \label{Fig: selected items for the naive semidefinite relaxation of (mKPC)}
   \end{minipage}
\end{figure}

The formulation given \refequation{eq: equivalence entre (mKPC) et la formulation SDP avec la matrice par blocs} allows us to build a first semidefinite relaxation, by dropping the non-convex rank-one constraint, giving the following semidefinite relaxation:

\begin{equation}
\left[\begin{array}{rl}
\text{minimize} & c^\top\diag\left(\mathbf{X}\right)\\
\text{ subject to} & w^\top\diag\left(\mathbf{X}\right)\geq q\\
& \forall i,j\in[\![n]\!]\text{ with }j-i>\Delta,\quad\left\lfloor\dfrac{j-i-1}{\Delta}\right\rfloor\mathbf{X}_{ij}\leq\displaystyle\sum\limits_{k=i+1}^{j-1}\mathbf{X}_{kk}\\
& \begin{pmatrix}
   1 & \diag\left(\mathbf{X}\right)^\top\\
   \diag\left(\mathbf{X}\right) & \mathbf{X}
\end{pmatrix}\succeq 0.
\end{array}\right.\tag{$(\text{mKPC})_{\text{SDP}}$}\label{eq: naive semidefinite relaxation}
\end{equation}
However, we observed that this first naive semidefinite relaxation does not provide in general better results than the linear relaxation \hyperref[eq: original formulation of the mKPC - LP version]{$\textcolor{CouleurCite}{\left(\text{mKPC}\right)_{\text{LP}}}$}. Let us illustrate this on one instance. \cref{Fig: selected items with the LP (mKPC),,Fig: selected items for the naive semidefinite relaxation of (mKPC)} show the profile of the selected objects for the linear and the naive semidefinite relaxation respectively. The bar corresponding to object $i\in[\![n]\!]$ is filled in proportion to the value of $x^\ast_i$ (or $\mathbf{X}^\ast_{ii}$), the optimal solution returned by the model.

Let's take a closer look at the solutions. For example, looking at $\mathbf{X}^\ast$, the optimal solution of \refequation{eq: naive semidefinite relaxation} and $x^\ast$ the optimal solution for the linear relaxation, between indexes $122$ and $126$ (around the location of the high peak):
\begin{equation}
   \mathbf{X}_{122-126}^\ast=\begin{pNiceMatrix}
& 0.9726 & 0.9751 & 0.9773 & 0.9726 & 0.9726 &\\
& 0.9751 & 0.9821 & 0.9847 & 0.9821 & 0.9821 &\\
& 0.9773 & 0.9847 & 0.9912 & 0.9912 & 0.9912 &\\
& 0.9726 & 0.9821 & 0.9912 & 1 & 1 &\\
& 0.9726 & 0.9821 & 0.9912 & 1 & 1 &\\
\CodeAfter
\tikz \filldraw[color=blue, fill=blue!5, very thick] ([xshift=-6.5mm,yshift=-2mm]1-2) rectangle ([xshift=6.5mm,yshift=2mm]1-2);
\tikz \draw (1-2) node[anchor=center]{$0.9726$} ;
\tikz \filldraw[color=blue, fill=blue!5, very thick] ([xshift=-6.5mm,yshift=-2mm]2-3) rectangle ([xshift=6.5mm,yshift=2mm]2-3);
\tikz \draw (2-3) node[anchor=center]{$0.9821$} ;
\tikz \filldraw[color=blue, fill=blue!5, very thick] ([xshift=-6.5mm,yshift=-2mm]3-4) rectangle ([xshift=6.5mm,yshift=2mm]3-4) ;
\tikz \draw (3-4) node[anchor=center]{$0.9912$} ;
\tikz \filldraw[color=blue, fill=blue!5, very thick] ([xshift=-6.5mm,yshift=-2mm]4-5) rectangle ([xshift=6.5mm,yshift=2mm]4-5) ;
\tikz \draw (4-5) node[anchor=center]{$1$} ;
\tikz \filldraw[color=blue, fill=blue!5, very thick] ([xshift=-6.5mm,yshift=-2mm]5-6) rectangle ([xshift=6.5mm,yshift=2mm]5-6) ;
\tikz \draw (5-6) node[anchor=center]{$1$} ;
\end{pNiceMatrix}
\qquad\qquad x^\ast_{122-126}=\begin{pmatrix}
0.9726 \\
0.9821 \\
0.9912 \\
1 \\
1 \\
\end{pmatrix}.\tag{$\left(\theEquationCounter\right)$}\label{example: apparent equality between the optimal solution of the LP and the optimal solution of the naive semidefinite relaxation}
\end{equation}\addtocounter{EquationCounter}{1}\null
This apparent equality between $\diag\left(\mathbf{X}^\ast\right)$ and $x^\ast$ can in fact be observed for all the other coordinates. Moreover, it seems that the coefficients off-diagonal $\mathbf{X}^\ast_{ij}$ are, within numerical error, derived from the product $x^\ast_ix^\ast_j$: indeed, for indexes $i\in[\![n]\!]$ such that $x^\ast_i=1$, we observe that $\mathbf{X}^\ast_{ij}=x_j^\ast$ for all $j\in[\![n]\!]$; and whenever $x^\ast_i=0$, we have $\mathbf{X}^\ast_{ij}=0$ for any $j\in[\![n]\!]$. This phenomenon can be observed in a very large number of instances and leads us to formulate the following property:

\begin{Property}[ROAD]\label{conj: La conjecture}
We will say that an instance of \emph{\refequation{eq: original formulation of the mKPC}} has the \emph{Rank-One with Adjusted Diagonal (ROAD) property} if the matrix
\begin{equation}
\mathbf{X}={x^\ast}{x^\ast}^\top+\Diag\left({x^\ast}-{x^\ast}^2\right)=\begin{pNiceArray}{cccc}[margin]
x^\ast_1 &  & \Block{2-2}<\large>{x^\ast_ix^\ast_j} &\\
&  & & \\
\Block{2-2}<\large>{x^\ast_ix^\ast_j} & & & \\
& & & x_n^\ast
\CodeAfter \line[shorten=6pt]{1-1}{4-4}
\end{pNiceArray}\tag{(ROAD)}\label{eq: ROAD idendity}
\end{equation}
is feasible for the (naive) semidefinite relaxation \emph{\refequation{eq: naive semidefinite relaxation}}, where $x^\ast$ is an optimal solution of \emph{\hyperref[eq: original formulation of the mKPC - LP version]{$\textcolor{CouleurCite}{\left(\text{mKPC}\right)_{\text{LP}}}$}} and where $\Diag:\R^n\rightarrow\SDP^n$ is the operator associating to a vector $u\in\R^n$ a matrix $\mathbf{U}\in\SDP^n$ with $u$ on the diagonal, and zero entries elsewhere.
\end{Property}

Empirical observations across a diverse range of instances suggest that most of them satisfy the \hyperref[conj: La conjecture]{ROAD property}. However, we can construct counterexamples where there is no clear relationship between the optimal solution of the linear relaxation and that of the semidefinite relaxation (see \ref{subsubsec: counterexample}).

Furthermore, even when considering an instance that satisfies the \hyperref[conj: La conjecture]{ROAD property}, the matrix $\mathbf{X}=x^\ast{x^\ast}^\top+\Diag(x^\ast-{x^\ast}^2)$ is feasible but not necessarily optimal for \refequation{eq: naive semidefinite relaxation}. For instance, in the example presented in \refequation{example: apparent equality between the optimal solution of the LP and the optimal solution of the naive semidefinite relaxation}, the optimal solution $\mathbf{X}^\ast$ of \refequation{eq: naive semidefinite relaxation} does not satisfy \refequation{eq: ROAD idendity}. Nevertheless, if we focus only on the objective values of the problems, we observe the following property that we leave as a conjecture:

\begin{Conjecture}\label{conj: conjecture 2 sur les bornes des problèmes}
$\Opt\left(\text{\emph{\refequation{eq: naive semidefinite relaxation}}}\right)\leq \Opt\left(\text{\emph{\hyperref[eq: original formulation of the mKPC - LP version]{\textcolor{CouleurCite}{$\left(\text{mKPC}\right)_{\text{LP}}$}}}}\right)\leq\Opt\emph{\refequation{eq: original formulation of the mKPC}}.$
\end{Conjecture}

In particular, for every instance that has the \hyperref[conj: La conjecture]{ROAD property}, we can use the optimal solution $x^\ast$ of the linear relaxation of \refequation{eq: original formulation of the mKPC} to build a solution $\mathbf{X}$ with the same objective value; hence, the optimal solution of the naive semidefinite relaxation has a lower objective value. Thus, \hyperref[conj: La conjecture]{ROAD property} implies \textbf{Conjecture~\ref{conj: conjecture 2 sur les bornes des problèmes}}. However, we do not know if \textbf{Conjecture \ref{conj: conjecture 2 sur les bornes des problèmes}} holds in general.

\subsection{Penalized SDP relaxation}

The idea of a penalty-based approach is to be able to construct good solutions by optimizing over a "score" of a solution that takes into account both accuracy and compactness, instead of minimizing the sum of the costs. In contrast, our approach introduces a penalty mechanism based on pairs of elements: it rewards configurations where the space between two indexes of a selected pair $(i,j)\in S^2$ is sufficiently filled by intermediate items $i<k<j$ while penalizing those that do not satisfy this criterion. 

\begin{Remark}
Other penalty-based approaches have been proposed, notably in \emph{\citeequation{curebal2024fixed}}, where a penalty is applied directly to the objective function when certain items are not selected.
\end{Remark}

Hence, we introduce a parameter $\lambda\in\R_+$ and change our objective function to deal with a balance between compactness and accuracy, and we come up with the following model:

\begin{equation}
   \left[\begin{array}{rl}
       \text{minimize} & c^\top\diag\left(\mathbf{X}\right)+\displaystyle\sum\limits_{\substack{1\leq i\leq j\leq n\\ j-i>\Delta}}\lambda\left(\left\lfloor\dfrac{j-i-1}{\Delta}\right\rfloor\mathbf{X}_{ij}-\displaystyle\sum\limits_{k=i+1}^{j-1}\mathbf{X}_{kk}\right) \\
       \text{subject to} & w^\top\diag\left(\mathbf{X}\right)\geq q\\
        & \begin{pmatrix}
           1 & \diag(\mathbf{X})^\top\\
           \diag(\mathbf{X}) & \mathbf{X}
       \end{pmatrix}\succeq 0.
   \end{array}\right.\tag{$\left(\text{PmKPC}\right)$}\label{eq: penality-based SDP}
\end{equation}

\begin{Remarks}[Alternative penalization approaches]
\item The model \emph{\refequation{eq: penality-based SDP}} does not correspond to the standard Lagrangian penalization and instead relies on an approach that rewards situations where the compactness constraints are verified. Indeed, if we consider the model obtained from building the Lagrangian penalization of \emph{\refequation{eq: naive semidefinite relaxation}}, we would get
$$\varphi_{\lambda}:\mathbf{X}\longmapsto c^\top\diag\left(\mathbf{X}\right)+\sum\limits_{\substack{1\leq i\leq j\leq n\\ j-i>\Delta}}\lambda\max\left(0,\left\lfloor\dfrac{j-i-1}{\Delta}\right\rfloor\mathbf{X}_{ij}-\sum\limits_{k=i+1}^{j-1}\mathbf{X}_{kk}\right).$$
Any $x\in[0,1]^n$ checking the knapsack constraint $w^\top x\geq q$ would yield a solution of the Lagrangian penalization model: applying the arguments from the beginning of \ref{subsubsec: counterexample}, the semidefinite positive matrix $\mathbf{X}:=xx^\top+\Diag(x-x^2)$ would check the knapsack constraint and the conic constraint. In particular, for the majority of instances that check the \emph{\hyperref[conj: La conjecture]{ROAD property}}, writing $x^\ast$ for the solution of \emph{\hyperref[eq: original formulation of the mKPC - LP version]{$\textcolor{CouleurCite}{\left(\text{mKPC}\right)_{\text{LP}}}$}}, we get that $\mathbf{X}:={x^\ast}{x^\ast}^\top+\Diag({x^\ast}-{x^\ast}^2)$ is a solution of the Lagrangian penalization model with the same objective value $\varphi_{\lambda}\left(\mathbf{X}\right)=c^\top {x^\ast}$. This is why we propose a model where verified compactness inequalities lead to a reward in the objective function thus avoids the problematic situations we encountered in the case of the naive semidefinite relaxation.
\item Another possible approach to penalization is to introduce individual penalty parameters $\lambda_{ij}$ for every pair $i,j\in [\![n]\!]$ in \emph{\refequation{eq: penality-based SDP}}. This extended parameterization could offer even more flexibility in tuning the objective function, allowing for more adjustments in response to the specific properties of the instances and the particular requirements of the problem. Here we prefer to stick with the simplest case with only one hyperparameter, to single out its role in the quality of obtained solutions (\emph{c.f.} \crefrange{Fig: spectrum: -log(lambda)=1}{Fig: spectrum: -log(lambda)=6} in \cref{sec: computationnal results}).
\end{Remarks}

To compare solution produced by different models (MIP, LP, SDP, penalized), we introduce metrics that assess solution quality independently of the objective function. These forthcoming metrics \refequation{def: parsimony} to \refequation{def: fractionnality} capture the key properties required to address the original statistical problem, namely accuracy and compactness. Note that the normalization included in \refequation{def: parsimony}, \refequation{def: compactness}, and \refequation{def: fractionnality} ensures that all three indicators take values in the range $[0,1]$, allowing for comparisons across instances of varying sizes.

Let $x^\ast\in[0,1]^n$ be a solution vector obtained from any of the models under consideration for a given instance. In particular, by \textbf{Proposition \ref{prop: formulation equivalente avec la contrainte conique}}, the solution vector produced by a semidefinite model can be directly extracted from the diagonal of the optimal solution matrix, \emph{i.e.}, $x^\ast=\diag\left(\mathbf{X}^\ast\right)$. To determine the set of selected items, we apply a simple rounding rule: an item $i$ is considered selected if and only if $x^\ast_i \geq 1/2$.

We then introduce the imprecision indicator as    
\begin{equation}
\mathbf{imp}\left(x^\ast\right):=\dfrac{c^\top x^\ast}{c^\top\mathbf{1}}\tag{$\left(\theEquationCounter\right)$}\label{def: parsimony}
\end{equation}\addtocounter{EquationCounter}{1}\null
where $c$ denotes the cost vector of the instance. This indicator corresponds to a normalized version of the original objective function of \refequation{eq: original formulation of the mKPC}: the aim of the original problem is to minimize the sum of the costs, which results in minimizing the imprecision of the solution.

Next, we consider a compactness indicator to evaluate the spatial distribution of the selected items: let $S\subseteq[\![n]\!]$ be the set of selected items for a given instance; the compactness indicator is defined as  
\begin{equation}
\mathbf{comp}\left(S\right):=\dfrac{1}{n}\max_{i,j\in S}\left\lbrace j-i-1\,\vert\,i,j \text{ consecutive}\right\rbrace\tag{$\left(\theEquationCounter\right)$}\label{def: compactness}
\end{equation}\addtocounter{EquationCounter}{1}\null
where $n$ denotes the total number of items in the instance. Finally, to quantify this ambiguity and evaluate how close a solution is to an integer vector, we introduce an additional metric: let $x^\ast$ be the solution vector of a model, the fractionality indicator is then defined as the normalized distance between $x^\ast$ and its closest $\left\lbrace 0,1\right\rbrace$-vector
\begin{equation}
\mathbf{frac}\left(x^\ast\right):=\dfrac{2}{\sqrt{n}}\left\Vert x^\ast-\left\lfloor x^\ast+\dfrac{1}{2}\cdot\mathbf{1}\right\rfloor\right\Vert_2\tag{$\left(\theEquationCounter\right)$}\label{def: fractionnality}
\end{equation}\addtocounter{EquationCounter}{1}\null
where the floor function is applied coordinate-wise to $x^\ast + \frac{1}{2} \mathbf{1}$. Since the considered formulations relax the integrality constraints on the solution vector, this fractionality indicator provides a relevant measure of solution quality and can be used to compare different models. Furthermore, as described in \citeequation{santini2024min}, \refequation{eq: original formulation of the mKPC} has instances that tend to produce relaxed solutions with many variables close to $\frac{1}{2}$, resulting in more branching when exploring the branch~{\&}~bound tree, making \refequation{def: fractionnality} even more relevant to evaluate the fractionnality of a solution, since $\mathbf{frac}(x)=1$ if and only if $x=\frac{1}{2} \mathbf{1}$.

Having a semidefinite formulation of the min-knapsack problem with compactness constraints and proposed a penalized version, we now turn our attention to several techniques to strenghen our SDP relaxations. In the following section, we introduce a series of valid inequalities and separation techniques aimed at reinforcing the product structure inherent to the binary formulation. These developments lead to tighter bounds and improved computational performance, as detailed in the subsequent numerical experiments.

\section{Strengthening the semidefinite relaxation}\label{sec: Strenghtening the semidefinite relaxation}

As discussed in \cref{subsec: a first sdp relaxation}, \refequation{eq: naive semidefinite relaxation} often fails to provide tighter bounds than its linear counterpart. In this section, we therefore introduce a series of valid inequalities, drawn from both semidefinite programming theory and combinatorial optimization, to strengthen the relaxation and enhance bound quality. We will then investigate how to separate a given fractional solution using the knapsack structure, and provide the theoretical justification for our approach.

\subsection{Valid inequalities}\label{subsec:Valid Inequalities}

We first recall that the proposed semidefinite reformulation \refequation{eq: equivalence entre (mKPC) et la formulation SDP avec la matrice par blocs} relies on the linearization of all quadratic terms from the original problem \refequation{eq: original formulation of the mKPC} by introducing a variable $\mathbf{X}_{ij}$ for each product $x_i x_j$, where $i,j$ range over $[\![n]\!]$. However, in the relaxed versions \refequation{eq: naive semidefinite relaxation} and \refequation{eq: penality-based SDP}, where the non-convex rank-one constraint has been dropped, the variables $\mathbf{X}_{ij}$ may no longer behave as true product terms.  

Inspired by techniques such as the Sherali-Adams cuts \citeequation{mathieu2009sherali}, we introduce valid inequalities designed to reinforce this product structure. These inequalities can be seamlessly integrated into both the standard semidefinite model and the penalty-based formulation, ensuring their applicability across different modeling approaches.

Let $I,J$ be disjoint subsets of $[\![n]\!]$ with $\vert I\cup J\vert=2$. Starting from the original $\left\lbrace 0,1\right\rbrace$ formulation of the \refequation{eq: original formulation of the mKPC}, we have that all the possible following inequalities
$$\prod_{i\in I}x_i\cdot\prod_{j\in J}\left(1-x_j\right)\geq 0$$
are quadratic valid inequalities for the  binary setup. Since $x\in\left\lbrace 0,1\right\rbrace^n$, and given the substitution $\mathbf{X} = xx^\top$ used to construct the semidefinite formulation, we can identify the product terms as follows: for any $i,j\in [\![n]\!]$, we associate $x_i^2$ with $\mathbf{X}_{ii}$ (\emph{resp.} $x_j^2$ with $\mathbf{X}_{jj}$) and $x_ix_j$ with $\mathbf{X}_{ij}$. This transformation lifts the inequalities into the matrix space, yielding the following set of valid inequalities for the semidefinite model:  
\begin{align}
\mathbf{X}_{ij}&\geq 0  &   \mathbf{X}_{ii}&\geq\mathbf{X}_{ij} &   \mathbf{X}_{ij}&\geq\mathbf{X}_{ii}+\mathbf{X}_{jj}-1\tag{$\left(\theEquationCounter\right)$}\label{eq: inegalites quadratiques 1}
\end{align}\addtocounter{EquationCounter}{1}\null
for all $i,j\in [\![n]\!]$. Moreover, for all $i,j,k\in[\![n]\!]$, the binary structure implies the following inequalities $(x_k-x_i)(x_k-x_j)\geq 0$ and $(1-x_k-x_i)(1-x_k-x_j)\geq 0$ which gives through the same lifting process, additional valid inequalities for the semidefinite setting:  
\begin{align}
\mathbf{X}_{kk}+\mathbf{X}_{ij}&\geq \mathbf{X}_{ik}+\mathbf{X}_{jk}  &   \mathbf{X}_{ik}+\mathbf{X}_{jk}+\mathbf{X}_{ij}&\geq\mathbf{X}_{ii}+\mathbf{X}_{jj}+\mathbf{X}_{kk}-1.\tag{$\left(\theEquationCounter\right)$}\label{eq: inegalites quadratiques 2}
\end{align}\addtocounter{EquationCounter}{1}\null
Furthermore, for any $j\in[\![n]\!]$ multiplying the knapsack constraint $w^\top x\geq q$ by $x_j$ and $1-x_j$ gives
\begin{align}
\sum_{i=1}^n w_i\mathbf{X}_{ij}&\geq q\mathbf{X}_{jj}  &   q\left(\mathbf{X}_{jj}-1\right)+\sum_{i=1}^nw_i\mathbf{X}_{ii}&\geq\sum_{i=1}^n w_i\mathbf{X}_{ij}.\tag{$\left(\theEquationCounter\right)$}\label{eq: inegalites quadratiques 3}
\end{align}\addtocounter{EquationCounter}{1}\null
Moreover, linerizing the Cauchy-Schwarz inequality $\left\vert\tr\left(\Diag(w)\mathbf{X}\right)\right\vert^2\leq\left\Vert\Diag\left(w\right)\right\Vert^2\left\Vert\mathbf{X}\right\Vert^2$ yields another quadratic valid inequality:
\begin{align}
\left(\sum_{i=1}^n w^2_i\right)\left(\sum_{1\leq i,j\leq n}\mathbf{X}_{ij}\right)&\geq \sum^n_{i=1}w_i^2\mathbf{X}_{ii}+2\sum_{1\leq i<j\leq n}w_iw_j\mathbf{X}_{ij}.\tag{$\left(\theEquationCounter\right)$}\label{eq: inegalites quadratiques 4}
\end{align}\addtocounter{EquationCounter}{1}\null

Overall, by incorporating inequalities \refequation{eq: inegalites quadratiques 1} to \refequation{eq: inegalites quadratiques 4}, to our semidefinite models, we get the strenghtened semidefinite relaxation $\textcolor{CouleurCite}{(\text{mKPC})^{+}_{\text{SDP}}}\label{eq: strenghened semidefinite relaxation}$ and its penalized version $\textcolor{CouleurCite}{\left(\text{PmKPC}\right)^+}\label{eq: strenghened penality-based SDP}$.


\subsection{Maximal insufficient subset cuts}

To strengthen again our semidefinite model, we investigate a class of valid inequalities derived from the classical \emph{cover inequalities} commonly used in knapsack problems. In our min-knapsack setting, the standard definition of cover should be inverted, motivating the definition of insufficient subsets: $S\subseteq [\![n]\!]$ is said to be \emph{insufficient} if $\sum_{i\in S}w_i<q$, that is to say that $S$ is not a valid selection of items for the knapsack constraint. 

In knapsack litterature, a cover set is a minimal set of items whose total weight exceeds the capacity, requiring at least one item to be removed. In our min-knapsack setting, this logic is reversed: we seek the largest possible subsets that remain insufficient, ensuring that at least one additional item must be included to satisfy the constraint. This motivates the following definition:

\begin{Definition}[Maximal insufficient subset]
If $S\subseteq[\![n]\!]$ is an insufficient subset, we call $S$ \emph{maximal} if: 
$$\forall j\notin S,\qquad w_j+\sum_{i\in S}w_i\geq q.$$
In words, $S$ is an insufficient subset that can be completed by any other item to verify the knapsack constraint. 
\end{Definition}

Now, if $S$ is such a maximal insufficient subset, we call \emph{maximal insufficient subset cut} \refequation{eq:Maximal Insufficient Subset Cut} associated to $S$ the following inequality: 
\begin{equation}  
\sum_{\substack{i\in[\![n]\!]\setminus S}}x_i\geq 1.\tag{(MISC)}\label{eq:Maximal Insufficient Subset Cut}  
\end{equation}
This inequality is valid for the classical knapsack problem. The bijection between the min-knapsack problem and the classical maximization knapsack \refequation{eq: Explicitation de la bijection entre min-knapsack et knapscak classique} also links the maximal insufficient subsets $S$ in the min-setting and the minimal covers $C$ of the complementary problem via $S\mapsto C:=[\![n]\!]\setminus S$. This fact allows us to recover all the theory from the classical knapsack theory and bring it into our setting. In particular, this tells us that any \refequation{eq:Maximal Insufficient Subset Cut} inequality is valid for \refequation{eq: mKP sans compactness} since, by definition, this inequality comes from a maximal insufficient subset $S$, meaning that at least one item outside $S$ must be selected to satisfy the knapsack constraint; and the \refequation{eq:Maximal Insufficient Subset Cut} formulation explicitly enforces this requirement.  

Furthermore, since this derivation relies only on the knapsack constraint itself, without assuming any specific structure in the distribution of selected items, \refequation{eq:Maximal Insufficient Subset Cut} remains valid for other mixed-integer models incorporating the constraint $w^\top x\geq q$. In particular, it applies to the semidefinite models we developed to address the change-point detection problem, including the "naive" relaxation \refequation{eq: naive semidefinite relaxation} and its strengthened version \hyperref[eq: strenghened semidefinite relaxation]{$\textcolor{CouleurCite}{(\text{mKPC})^{+}_{\text{SDP}}}$}, as well as the penalty-based formulation \refequation{eq: penality-based SDP} and its strengthened counterpart \hyperref[eq: strenghened penality-based SDP]{$\textcolor{CouleurCite}{\left(\text{PmKPC}\right)^+}$}.

In particular, we will know focus on the semidefinite framework where the \refequation{eq:Maximal Insufficient Subset Cut} is formulated as:
\begin{equation}  
\sum_{\substack{i\in[\![n]\!]\setminus S}}\mathbf{X}_{ii}\geq 1.\tag{$(\text{MISC})_{\text{SDP}}$}\label{eq:Maximal Insufficient Subset Cut SDP}  
\end{equation}  
We can now show how to integrate those maximal insufficent subset cuts in a cutting plane scheme: given $\mathbf{X}^\ast$ a fractional solution of one of our semidefinite formulation, we want to find $S\subseteq[\![n]\!]$ a maximal insufficient subset that yield an inequality \refequation{eq:Maximal Insufficient Subset Cut SDP} that separates $\mathbf{X}^\ast$, or certify that no such subset exists.

To adress this problem, we note that $S\subseteq[\![n]\!]$ is a maximal insufficient subset that yields a cut \refequation{eq:Maximal Insufficient Subset Cut SDP} separating $\mathbf{X}^\ast$ if and only if $S$ is such that
$$\sum_{i\in [\![n]\!]\setminus S}\mathbf{X}^\ast_{ii}<1\qquad\qquad\text{or, equivalently,}\qquad\qquad\sum_{i=1}^n\left(1-\mathds{1}_{S}(i)\right)\mathbf{X}^\ast_{ii}<1.$$
Introducing binary variables $\alpha_i$ for $i\in[\![n]\!]$ with value $1$ if and only if $i\in S$, we build the following separation problem:
\begin{equation}
    \left[\begin{array}{rl}
        \text{minimize} & \displaystyle\sum\limits_{i=1}^n\left(1-\alpha_i\right) \mathbf{X}^\ast_{ii}\\
        \text{subject to} & \displaystyle\sum\limits^n_{i=1}w_i\alpha_i\leq q-\varepsilon\\
        & \alpha\in\lbrace 0,1\rbrace^n.
    \end{array}\right.\tag{$\left(\theEquationCounter\right)$}\label{eq: separation problem}
\end{equation}\addtocounter{EquationCounter}{1}\null
where $\varepsilon>0$ is small enough to make sure that the feasible solutions of the integer problems gives exactly the set of all insufficient subsets:
$$\left\lbrace \alpha\in\R^n\,\left\vert\,w^\top\alpha<q\right.\right\rbrace\cap\left\lbrace 0,1\right\rbrace^n=\left\lbrace \alpha\in\R^n\,\left\vert\,w^\top\alpha\leq q-\varepsilon\right.\right\rbrace\cap\left\lbrace 0,1\right\rbrace^n.$$

With this formulation, \refequation{eq: separation problem} is a combinatorial problem that aim to give a insufficient subset yielding a cut separating $\mathbf{X}^\ast$, or decide whether such set exists. Indeed, if $\Opt\text{\refequation{eq: separation problem}}\geq 1$, then we know that there is no insufficient subset, and in particular no maximal insufficient subset, that yield a separating \refequation{eq:Maximal Insufficient Subset Cut SDP} inequality. In the other hand, if $\Opt\text{\refequation{eq: separation problem}}<1$, the optimal solution $\alpha^\ast$ gives a set $S=\left\lbrace i\text{ such that }\alpha^\ast_i=1\right\rbrace$ such that $\mathbf{X}^\ast$ violates \refequation{eq:Maximal Insufficient Subset Cut SDP}. 

Furthermore, rewritting the objective of \refequation{eq: separation problem} as a maximization problem, one can observe that it is precisely equivalent to a classical knapsack problem. Given that such problems are often considered "easy" to solve optimally, the insufficient subset $S=\left\lbrace i\text{ such that }\alpha^\ast_i=1\right\rbrace$ can therefore be computed in pseudo-polynomial time \citeequation{pisinger2005hard}. For further details and techniques related to knapsack problems, see, \emph{e.g.}, \citeequation{hojny2020knapsack}.

However, the solution set $S=\left\lbrace i\text{ such that }\alpha^\ast_i=1\right\rbrace$ returned by \emph{\refequation{eq: separation problem}} is not necessarly maximal. The maximality of the subset could be ensured with additional "big $M$ constraints"
$$\forall j\in[\![n]\!],\qquad M\left(1-\alpha_j\right)+w_j+\sum^n_{i=1}w_i\alpha_i\geq q$$
but adding such constraints would bring a lot of complexity to the separation problem, meaning that we would lose our ability to compute $S=\left\lbrace i\text{ such that }\alpha^\ast_i=1\right\rbrace$ easely. Instead, we will "correct" the returned set using the following greedy algorithm:

\begin{algorithm}
\begin{algorithmic}
\State $\overline{S}\gets\left\lbrace i\text{ such that }\alpha^\ast_i=1\right\rbrace$ 
\While{$\displaystyle\sum\limits_{i\in \overline{S}} w_i< q$}
\State add to $\overline{S}$ the lightest item from $[\![n]\!]\setminus\overline{S}$\Comment{\emph{repeats until $\overline{S}$ is no more insufficient}}
\EndWhile
\State remove the last item added to $\overline{S}$
\State \textbf{return} $\overline{S}$
\end{algorithmic}
\caption{greedy approach to correct $S$}
\label{alg: algorithm to make a subset maximal}
\end{algorithm}

Our motivation to rebuild maximality out of the output of \refequation{eq: separation problem} comes from the fact that we want to have the tighter cut possible in order to leave behing the largest possible amount of fractionnal feasible points. Using the bijection between the min-knapsack problem and its complementary described earlier, we can state the following proposition:

\begin{Proposition}[Adapted from the classical knapsack litterature, see \emph{e.g.} \textbf{Proposition 7.1} in \citeequation{conforti2014integer}]
Let $S$ be an insufficient subset for \emph{\refequation{eq: mKP sans compactness}}, then the insufficient subset inequality
$$\sum_{\substack{i\in[\![n]\!]\setminus S}}x_i\geq 1$$
is facet-defining for $\conv\left(\mathcal{K}\right)\cap\left\lbrace x\in\R^n\,\vert\,\forall i\in S, x_i=1\right\rbrace$, where $\mathcal{K}$ denotes the set of all integer feasible points of \emph{\refequation{eq: mKP sans compactness}}.
\end{Proposition}

Hence, the maximality of an insufficient subset is a relevant property to ask for to get the tightest inequality possible. We now have the following result:

\begin{Proposition}\label{prop: propriétés des MIS}
Let $\alpha^\ast$ be an optimal solution of \refequation{eq: separation problem}, and $\overline{S}$ the corrected set returned by \emph{\textbf{Algorithm \ref{alg: algorithm to make a subset maximal}}}.
\begin{enumerate}[label=\emph{(\roman*)}]
\item The set $\overline{S}$ is a maximal insufficient subset.
\item If $\Opt\text{\refequation{eq: separation problem}}<1$, then the \emph{\refequation{eq:Maximal Insufficient Subset Cut SDP}} associated to $\overline{S}$ separates $\mathbf{X}^\ast$.
\item Otherwise, if $\Opt\text{\refequation{eq: separation problem}}\geq 1$, then there exists no maximal insufficient maximal subset $S$ such that \emph{\refequation{eq:Maximal Insufficient Subset Cut SDP}} associated to $\overline{S}$ separates $\mathbf{X}^\ast$.
\end{enumerate}
\end{Proposition}
\begin{proof}
The point (iii) has already been discussed after the presentation of the separation problem \refequation{eq: separation problem}. Let $S:=\left\lbrace i\text{ such that }\alpha^\ast_i=1\right\rbrace$ be the input set of \textbf{Algorithm \ref{alg: algorithm to make a subset maximal}}. First, notice that \textbf{Algorithm \ref{alg: algorithm to make a subset maximal}} runs in a finite number of steps since $w_1+\dots+w_n\geq q$: we therefore cannot add items to $S$ indefinitely in such a way that the resulting set $\overline{S}$ remains insufficient.
By construction, at each iteration of the while loop, an item is added to $\overline{S}$. Let $\mathbf{i}_k$ be the item we added at iteration $k$, and let $m$ be the total number of iterations. The returned set is then $\overline{S}=S\sqcup\left\lbrace\mathbf{i}_1,\dots,\mathbf{i}_{m-1}\right\rbrace$ and it is obviously insufficient, otherwise the loop would have stopped earlier. Moreover, by definition of $\mathbf{i}_m$, for each $j\in[\![n]\!]\setminus\overline{S}$, we have $w_j\geq w_{\mathbf{i}_m}$, therefore
$$w_j+\sum_{i\in\overline{S}}w_i\geq \sum_{i\in\overline{S}\sqcup\left\lbrace \mathbf{i}_m\right\rbrace}w_i\geq q$$
meaning that $\overline{S}$ is maximal, showing (i).

Now if $\Opt\text{\refequation{eq: separation problem}}<1$, then
$$\sum_{i\in [\![n]\!]\setminus\overline{S}}\mathbf{X}_{ii}^\ast=\sum_{\substack{i\in [\![n]\!]\setminus S\\ i\neq\mathbf{i}_1,\dots,\mathbf{i}_{m-1}}}\mathbf{X}_{ii}^\ast=\underbrace{\sum_{i=1}^n\left(1-\alpha_i^\ast\right)\mathbf{X}_{ii}^\ast}_{<1}-\underbrace{\sum_{k=1}^{m-1}\mathbf{X}_{\mathbf{i}_k\mathbf{i}_k}^\ast}_{\geq 0}<1$$
so the \refequation{eq:Maximal Insufficient Subset Cut SDP} associated to $\overline{S}$ separates $\mathbf{X}^\ast$, and the label (ii) is shown.
\end{proof}

We can see that the global procedure which constructs a maximal insufficient subset separating a fractional point $\mathbf{X}^\ast$ is pseudo-polynomial, and therefore takes a relatively small computation time compared to the global computation time of the various models. This allows us, for one of the proposed semidefinite models \hyperref[eq: strenghened semidefinite relaxation]{$\textcolor{CouleurCite}{(\text{mKPC})^{+}_{\text{SDP}}}$} or \hyperref[eq: strenghened penality-based SDP]{$\textcolor{CouleurCite}{\left(\text{PmKPC}\right)^+}$}, to get an optimal solution $\mathbf{X}^\ast$ and try to get a separating inequality that can be added to the model. Theorically, it is possible to have instances where no separating \refequation{eq:Maximal Insufficient Subset Cut} exists, that is to say that $\Opt\text{\refequation{eq: separation problem}}\geq 1$, hence we will integrate in our separation heuristic a solve of the linear relaxation of \refequation{eq: separation problem}: indeed, if the objective of the linear relaxation of \refequation{eq: separation problem} is greater or equal to $1$, we then know that the value of the integer problem is greater or equal to $1$, implying that no maximal insufficient subset separating $\mathbf{X}^\ast$ exists, by \textbf{Proposition \ref{prop: propriétés des MIS}}. 

We can now describe our separating procedure for a semidefinite model $M$, as proposed in \textbf{Algorithm \ref{alg: separation procedure}}.

\begin{algorithm}
\begin{algorithmic}
\State $\mathbf{X}^\ast\gets\argmin\left(M\right)$
\State solve the linear relaxation of \refequation{eq: separation problem}
\If{$\Opt\left(\text{\refequation{eq: separation problem}}_{\textcolor{CouleurCite}{\text{LP}}}\right)\geq 1$}\Comment{\emph{$\Opt\text{\refequation{eq: separation problem}}\geq\Opt\left(\text{\refequation{eq: separation problem}}_{\textcolor{CouleurCite}{\text{\emph{LP}}}}\right)\geq 1$, so $\not\exists$ \emph{\refequation{eq:Maximal Insufficient Subset Cut SDP}} separating $\mathbf{X}^\ast$}}
\State \textbf{return} $\mathbf{X}^\ast$
\Else
\State solve \refequation{eq: separation problem}
\If{$\Opt\text{\refequation{eq: separation problem}}\geq 1$}\Comment{\emph{$\not\exists$\emph{ \refequation{eq:Maximal Insufficient Subset Cut SDP}} separating $\mathbf{X}^\ast$}}
\State \textbf{return} $\mathbf{X}^\ast$
\Else
\State run \textbf{Algorithm \ref{alg: algorithm to make a subset maximal}} with imput $\left\lbrace i\text{ such that }\alpha^\ast_i=1\right\rbrace$
\State get $\overline{S}$\Comment{\emph{by \emph{\textbf{Proposition \ref{prop: propriétés des MIS}}}, $\overline{S}$ is maximal \& insuffient, such that \emph{\refequation{eq:Maximal Insufficient Subset Cut SDP}} separates $\mathbf{X}^\ast$}}
\State solve current model with the cut \refequation{eq:Maximal Insufficient Subset Cut SDP} 	associated to $\overline{S}$
\EndIf
\EndIf
\end{algorithmic}
\caption{separation procedure}
\label{alg: separation procedure}
\end{algorithm}


In this section, we have developped strengthening techniques for our semidefinite relaxations, valid inequalities, by lifting classical binary quadratic constraints to get valid inequalities, and by incorporating the combinatorial structure of the problem in a separation procedure. Now, we illustrate the effect of these inequalities on the quality of the bounds through numerical comparisons.

\section{Computationnal results}\label{sec: computationnal results}

This section presents the results of computational experiments conducted to evaluate the effectiveness of the proposed models.

\subsection{Experimental setup}\label{subsec: Instances}

The implementation was done in Julia 1.11.2 \citeequation{bezanson2017julia} with the optimization package JuMP.jl \citeequation{Lubin2023}, version 1.23.6; using HiGHS 1.13.0 \citeequation{HiGHSjl} as the MIP solver and Mosek 10.2.0 \citeequation{MOSEKjl} for semidefinite programming. All experiments were performed on Grid'5000 cluster machines \citeequation{cappello2005grid}, equipped with ${\text{Intel}}^{\text{\textregistered}}$ Xeon Gold 6130 CPU running at 2.10 GHz. A time limit of 7200 seconds was imposed on all solvers.

Since the computational difficulty is highly dependent on the instance characteristics \citeequation{santini2024min}, and many instances of \refequation{eq: original formulation of the mKPC} can be efficiently solved using the methods proposed in \citeequation{santini2024min}, our numerical experiments focused on the most challenging cases, see \ref{app: Generation of instances}.

\subsection{Computational experiments}

We begin by analyzing the performance of the unpenalized semidefinite relaxation \refequation{eq: naive semidefinite relaxation} and its strengthened variant \hyperref[eq: strenghened semidefinite relaxation]{$\textcolor{CouleurCite}{(\text{mKPC})^{+}_{\text{SDP}}}$}, comparing them to both the MIP formulation of \refequation{eq: original formulation of the mKPC} and its linear relaxation.

\begin{figure}[H]
	\begin{minipage}{0.48\linewidth}
		\centering
    		\includegraphics[width=\linewidth,keepaspectratio]{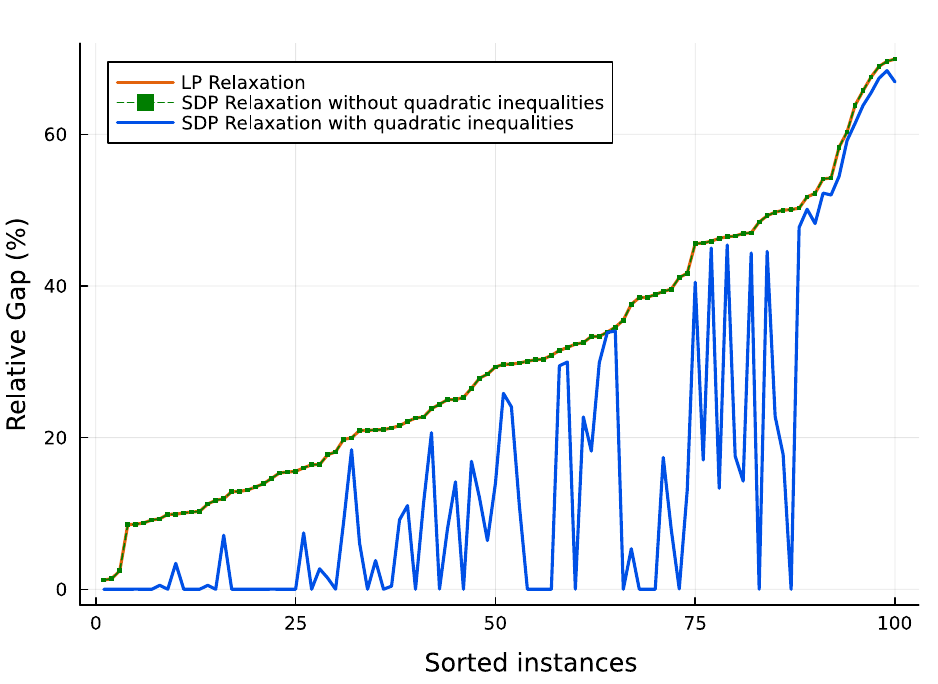}
   \caption{Relative gap between the MIP model \refequation{eq: original formulation of the mKPC} and its linear relaxation, the SDP relaxation \refequation{eq: naive semidefinite relaxation} and its strengthened version \hyperref[eq: strenghened semidefinite relaxation]{$\textcolor{CouleurCite}{(\text{mKPC})^{+}_{\text{SDP}}}$}.}\label{Fig: relative gap vs LP vs SDP sans et avec inegalités quadratiques}
   \end{minipage}\hfill
   \begin{minipage}{0.48\linewidth}
       \centering
   		\includegraphics[width=\linewidth,keepaspectratio]{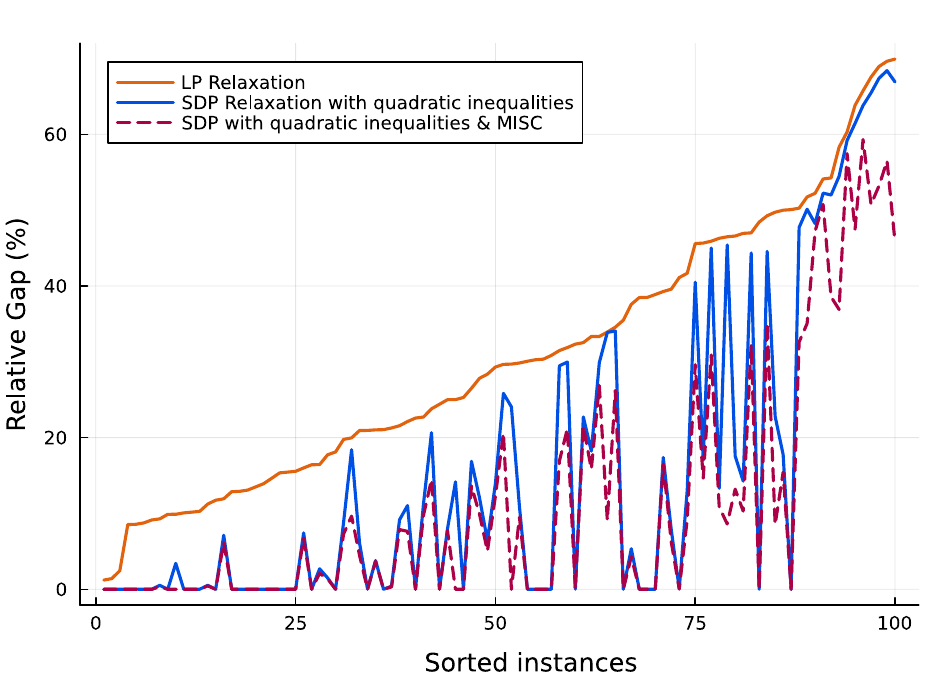}
   		\caption{Relative gap between the MIP model \refequation{eq: original formulation of the mKPC} and its linear relaxation, the strengthened SDP relaxation \hyperref[eq: strenghened semidefinite relaxation]{$\textcolor{CouleurCite}{(\text{mKPC})^{+}_{\text{SDP}}}$} and the SDP with a separating \refequation{eq:Maximal Insufficient Subset Cut SDP} using \textbf{Algorithm \ref{alg: separation procedure}}.}
   	\label{Fig: relative gap vs LP vs SDP avec et sans MISC}
   \end{minipage}
\end{figure}

\paragraph{Strenghening valid inequalities \refequation{eq: inegalites quadratiques 1} to \refequation{eq: inegalites quadratiques 4}} \cref{Fig: relative gap vs LP vs SDP sans et avec inegalités quadratiques} illustrates the relative gap, defined as
$$\text{\textsc{Gap}}\left(\%\right)=100\cdot\dfrac{\text{\textsc{ub}}-\text{\textsc{lb}}}{\text{\textsc{ub}}}$$
where \textsc{ub} corresponds to the objective value obtained by the MIP model, and \textsc{lb} denotes the bound provided by the studied relaxation. \cref{Fig: relative gap vs LP vs SDP sans et avec inegalités quadratiques} compares the relative gap of the MIP formulation \refequation{eq: original formulation of the mKPC} with its linear relaxation, as well as with the "naive’ semidefinite relaxation \refequation{eq: naive semidefinite relaxation} and its strengthened counterpart \hyperref[eq: strenghened semidefinite relaxation]{$\textcolor{CouleurCite}{(\text{mKPC})^{+}_{\text{SDP}}}$}, over a benchmark of 100 instances. The instances are ordered by increasing relative gap in the linear relaxation, providing an approximate classification of computational difficulty. The near-perfect overlap between the linear relaxation curve (\textcolor{CouleurLP}{\textbf{-----}}) and the "naive’ semidefinite relaxation curve (\textcolor{CouleurSDPNaif}{\textbf{-~\scriptsize{\faSquare}\normalsize{~-}}}) empirically confirms \textbf{Conjecture \ref{conj: conjecture 2 sur les bornes des problèmes}} from Section~\ref{subsec: a first sdp relaxation}. More importantly, \cref{Fig: relative gap vs LP vs SDP sans et avec inegalités quadratiques} underscores the impact of the quadratic inequalities introduced in the strengthened model \hyperref[eq: strenghened semidefinite relaxation]{$\textcolor{CouleurCite}{(\text{mKPC})^{+}_{\text{SDP}}}$}. The corresponding curve (\textcolor{CouleurSDPQuadraticInequalities}{\textbf{-----}}) consistently lies below the linear relaxation curve, demonstrating that incorporating quadratic inequalities \refequation{eq: inegalites quadratiques 1} to \refequation{eq: inegalites quadratiques 4} effectively strengthens the bound, even for the most challenging instances located on the right-hand side of the graph. Additionally, the strengthened semidefinite model frequently achieves a relative gap of $0$, indicating that the instance has been solved to optimality.

\paragraph{Separating \emph{\refequation{eq:Maximal Insufficient Subset Cut}} strategy}  
\cref{Fig: relative gap vs LP vs SDP avec et sans MISC} illustrates the relative gap to the MIP model for the linear relaxation, the strengthened semidefinite formulation \hyperref[eq: strenghened semidefinite relaxation]{$\textcolor{CouleurCite}{(\text{mKPC})^{+}_{\text{SDP}}}$}, and the bound obtained using the separation procedure \textbf{Algorithm \ref{alg: separation procedure}}. We observe that incorporating a separating \refequation{eq:Maximal Insufficient Subset Cut SDP} strategy into the strengthened semidefinite relaxation via \textbf{Algorithm \ref{alg: separation procedure}} leads to tighter objective bounds (\textcolor{CouleurSDPMISC}{\textbf{-~-~-}}). Specifically, for instances where the strengthened semidefinite model does not achieve a zero relative gap, \textbf{Proposition \ref{prop: formulation equivalente avec la contrainte conique}} implies that the returned optimal solution $\mathbf{X}^\ast$ contains fractional coefficients, and in such cases, the separation problem \refequation{eq: separation problem} aims to provide a valid inequality that cuts $\mathbf{X}^\ast$, further improving the bound.

\bigskip

\begin{figure}[H]
\centering
   \begin{minipage}{0.58\linewidth}
       \centering
       \includegraphics[width=\linewidth,keepaspectratio]{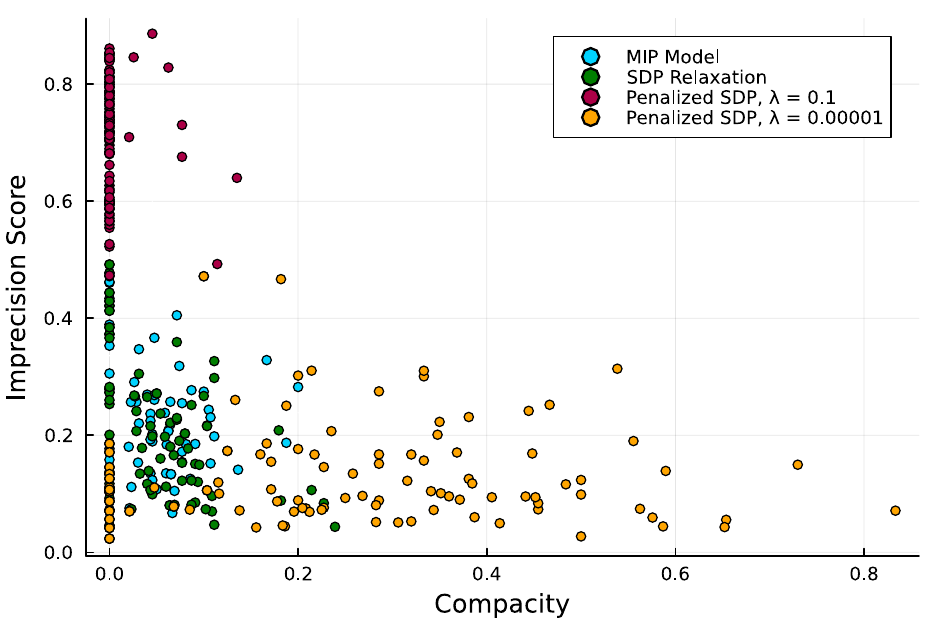}
       \vspace*{-8mm}\caption{plot of the balance between the imprecision score and compacity for the MIP model, the strengthened semidefinite relaxation \hyperref[eq: strenghened semidefinite relaxation]{$\textcolor{CouleurCite}{(\text{mKPC})^{+}_{\text{SDP}}}$} and two penalized \hyperref[eq: strenghened penality-based SDP]{$\textcolor{CouleurCite}{\left(\text{PmKPC}\right)^+}$}.}
       \label{Fig: spectrum: compacity vs parsimony on the benchmark}
   \end{minipage}\hfill
   \begin{minipage}{0.38\linewidth}
       \centering
       \begin{longtable}[c]{|cc|c|}
       \hline
       \multicolumn{2}{|c|}{\textsc{Model}} &
         \begin{tabular}[c]{@{}c@{}}\textsc{Average}\\ \textsc{fractionnality}
         \end{tabular} \\ \hline\hline
       \endhead
       \hline
       \endfoot
       \endlastfoot
       \multirow{4}{*}{\rotatebox{90}{\hyperref[eq: strenghened penality-based SDP]{$\textcolor{CouleurCite}{\left(\text{PmKPC}\right)^+}$}}} &
         $\lambda=\exposantvanilla{0}$ & $2.779\exposant{-3}$\Tstrut \\
        & $\lambda=\exposantvanilla{-2}$ & $2.339\exposant{-3}$ \\
        & $\lambda=\exposantvanilla{-4}$ & $1.673\exposant{-2}$ \\
        & $\lambda=\exposantvanilla{-6}$ & $4.326\exposant{-2}$ \\
        \hline
       \multicolumn{2}{|c|}{LP}    & $2.583\exposant{-1}$\Tstrut \\
       \multicolumn{2}{|c|}{\hyperref[eq: strenghened semidefinite relaxation]{$\textcolor{CouleurCite}{(\text{mKPC})^{+}_{\text{SDP}}}$}} &
         $1.904\exposant{-1}$ \\[1ex]
         \hline
       \end{longtable}
       \vspace{15pt}
       \caption{average fractionnality observed on the benchmark for (LP), \hyperref[eq: strenghened semidefinite relaxation]{$\textcolor{CouleurCite}{(\text{mKPC})^{+}_{\text{SDP}}}$} and \hyperref[eq: strenghened penality-based SDP]{$\textcolor{CouleurCite}{\left(\text{PmKPC}\right)^+}$} for different $\lambda$.}
       \label{Fig: tab: average fractionnality on the benchmark}
   \end{minipage}
\end{figure}

\begin{figure}[H]
\centering
\begin{subcaptiongroup}
       \begin{minipage}{0.33\linewidth}
       \centering
       \includegraphics[width=\linewidth,keepaspectratio]{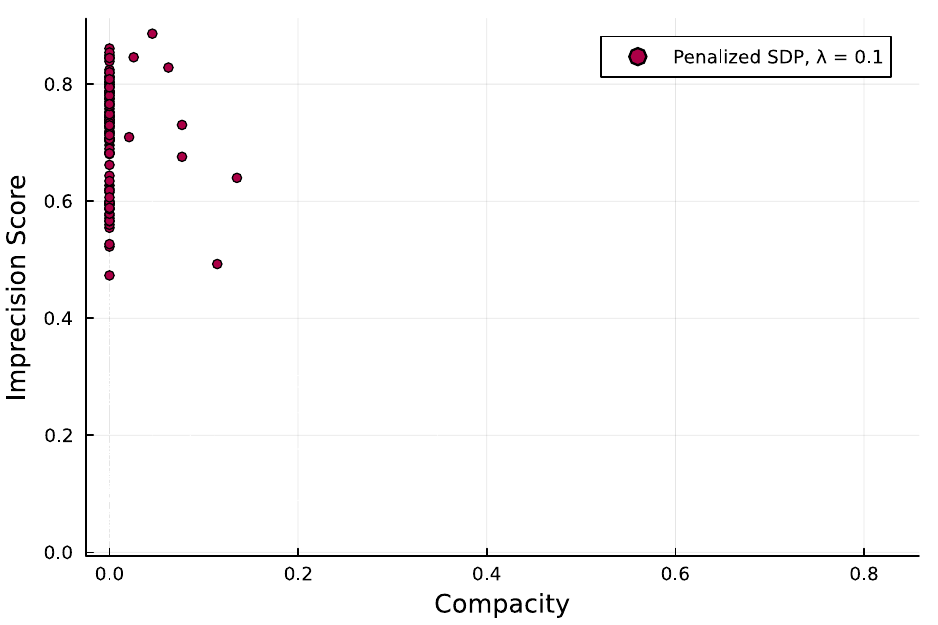}
       \caption{$\lambda=\exposantvanilla{-1}$}
       \label{Fig: spectrum: -log(lambda)=1}
   \end{minipage}\hfill
       \begin{minipage}{0.33\linewidth}
       \centering
       \includegraphics[width=\linewidth,keepaspectratio]{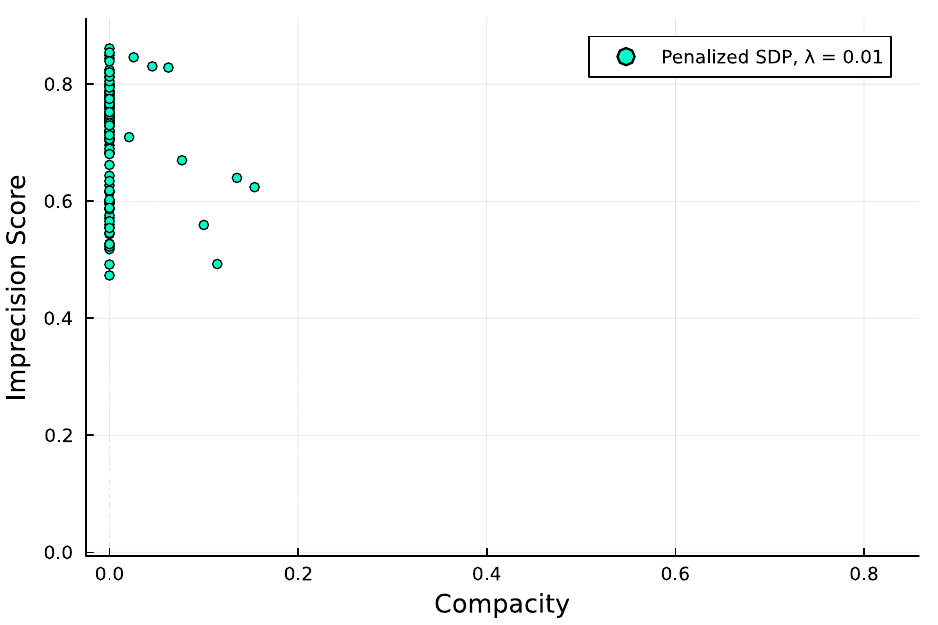}
       \caption{$\lambda=\exposantvanilla{-2}$}
       \label{Fig: spectrum: -log(lambda)=2}
   \end{minipage}
       \begin{minipage}{0.33\linewidth}
       \centering
       \includegraphics[width=\linewidth,keepaspectratio]{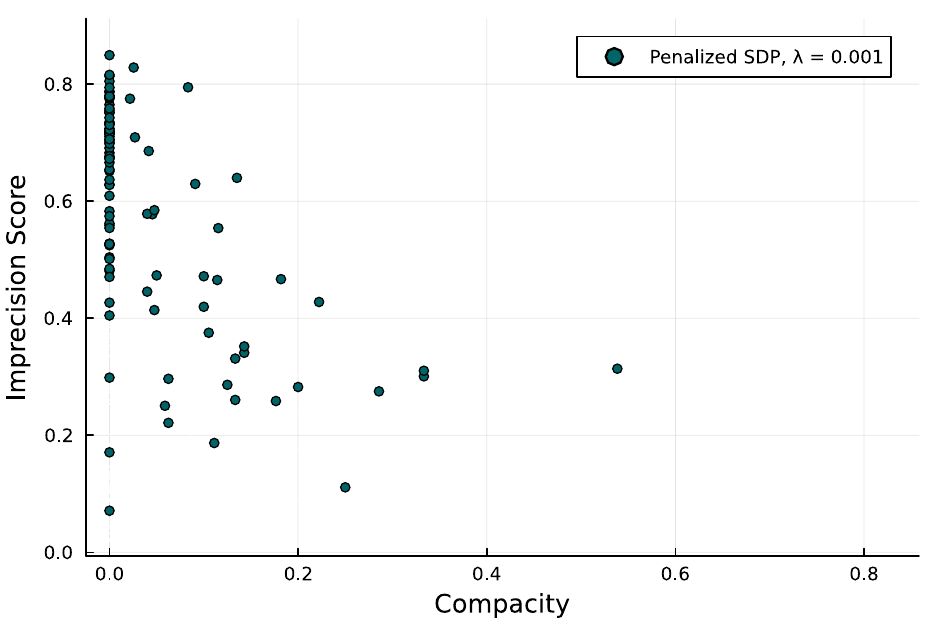}
       \caption{$\lambda=\exposantvanilla{-3}$}
       \label{Fig: spectrum: -log(lambda)=3}
   \end{minipage}\\
\vspace*{0.5cm}
       \begin{minipage}{0.33\linewidth}
       \centering
       \includegraphics[width=\linewidth,keepaspectratio]{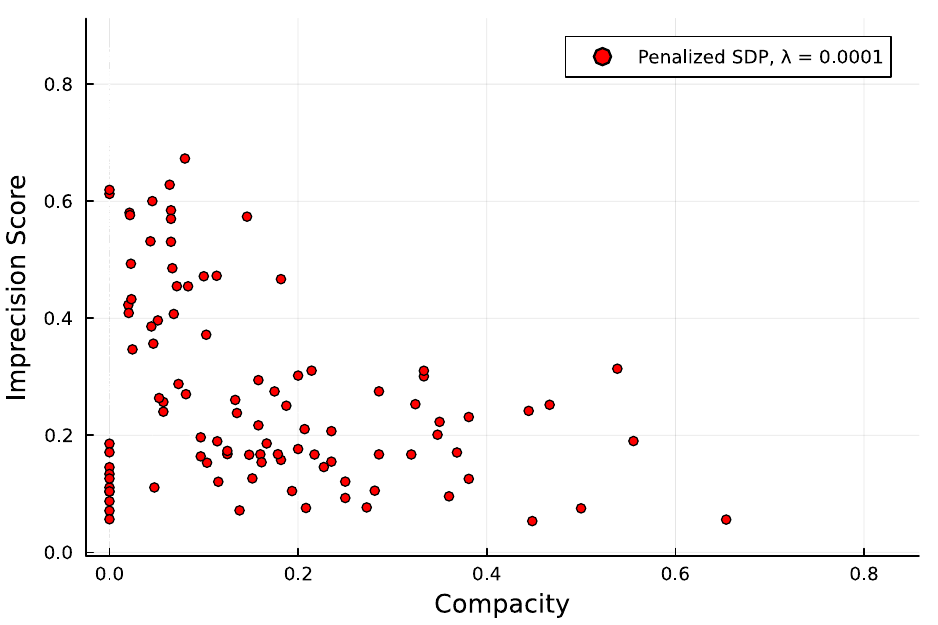}
       \caption{$\lambda=\exposantvanilla{-4}$}
       \label{Fig: spectrum: -log(lambda)=4}
   \end{minipage}\hfill
       \begin{minipage}{0.33\linewidth}
       \centering
       \includegraphics[width=\linewidth,keepaspectratio]{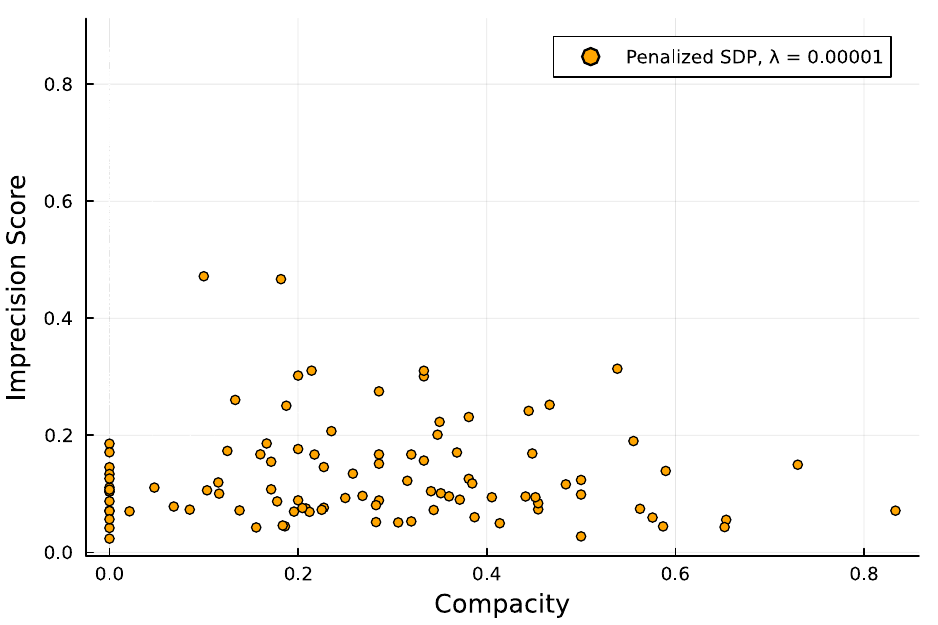}
       \caption{$\lambda=\exposantvanilla{-5}$}
       \label{Fig: spectrum: -log(lambda)=5}
   \end{minipage}\hfill
       \begin{minipage}{0.33\linewidth}
       \centering
       \includegraphics[width=\linewidth,keepaspectratio]{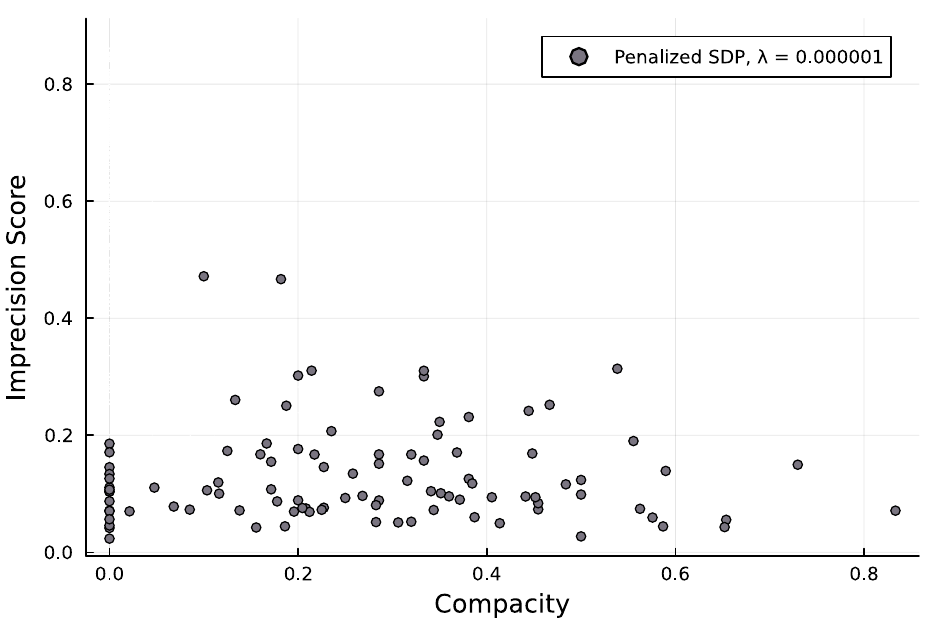}
       \caption{$\lambda=\exposantvanilla{-6}$}
       \label{Fig: spectrum: -log(lambda)=6}
   \end{minipage}
   \end{subcaptiongroup}
   \caption{Plots for different values of $\lambda$}
\end{figure}

\paragraph{Penalized semidefinite approach} When analyzing penalized versions of the model, it is essential to consider the metrics of imprecision, compactness, and fractionality, as introduced in \refequation{def: parsimony}, \refequation{def: compactness} and \refequation{def: fractionnality}. \cref{Fig: spectrum: compacity vs parsimony on the benchmark} illustrates the behavior of the studied models by representing each solved instance as a point whose coordinates correspond to the compactness and imprecision scores of the obtained solution: each set of dots of the same color stands for a studied model. In order to avoid poor visibility on the graph, we only show the data for the strenghened penalized version \hyperref[eq: strenghened penality-based SDP]{$\textcolor{CouleurCite}{\left(\text{PmKPC}\right)^+}$} for two extreme values of $\lambda$ which would give rise to characteristic behavior. This visualization highlights the influence of different values of $\lambda$ on the trade-off between accuracy and compactness: the set of dots corresponding to \refequation{eq: penality-based SDP} with $\lambda=\exposantvanilla{-1}$ on \cref{Fig: spectrum: compacity vs parsimony on the benchmark} indicates that setting $\lambda=\exposantvanilla{-1}$ tends to produce highly compact solutions, but at the expense of an increasing imprecision score, indicating that this value overly penalizes compactness. Conversely, $\lambda=\exposantvanilla{-6}$ favors accuracy but results in insufficient compactness, suggesting that the penalty is too weak. The plots in \cref{Fig: spectrum: compacity vs parsimony on the benchmark} provide an overview of these effects, while a more detailed impact of varying $\lambda$ values, for $\lambda\in \left\lbrace \exposantvanilla{-1}, \exposantvanilla{-2},\dots, \exposantvanilla{-6} \right\rbrace$, can be observed in \crefrange{Fig: spectrum: -log(lambda)=1}{Fig: spectrum: -log(lambda)=6}. For example, with this instance size, a balanced trade-off between accuracy and compactness can be achieved with $\lambda\approx\exposantvanilla{-3}$, as illustrated in \cref{Fig: spectrum: -log(lambda)=3}.

\paragraph{Fractionality}  
From a fractionality perspective, the table in \cref{Fig: tab: average fractionnality on the benchmark} reports the average fractionality across the entire benchmark for the different models considered. The results indicate that the penalized models \hyperref[eq: strenghened penality-based SDP]{$\textcolor{CouleurCite}{\left(\text{PmKPC}\right)^+}$} tend to yield solutions that are closer to integer vectors than those obtained from the relaxation \hyperref[eq: strenghened semidefinite relaxation]{$\textcolor{CouleurCite}{(\text{mKPC})^{+}_{\text{SDP}}}$}. Furthermore, within the penalized models, increasing the value of $\lambda$ generally leads to solutions with lower fractionality. However, even after selecting $\lambda$ to strike a balance between accuracy and compactness, fractionality can still be further reduced by incorporating a \refequation{eq:Maximal Insufficient Subset Cut} separating inequality strategy. \cref{Fig: fractionnality with the separation procedure} illustrates the average fractionality across the benchmark for the penalized model \hyperref[eq: strenghened penality-based SDP]{$\textcolor{CouleurCite}{\left(\text{PmKPC}\right)^+}$} with $\lambda\in\left\lbrace\exposantvanilla{0},\exposantvanilla{-2},\exposantvanilla{-4},\exposantvanilla{-6}\right\rbrace$, compared to the fractionality observed when a single separating \refequation{eq:Maximal Insufficient Subset Cut SDP} cut is added to the model. Notably, since the considered instances tend to produce solutions whose coordinates are close to $\frac{1}{2}$ \citeequation{santini2024min}, leading to a fractionality score $\mathbf{frac}\left(\diag\left(\mathbf{X}^\ast\right)\right)$ close to $1$, \cref{Fig: fractionnality with the separation procedure} demonstrates that combining the penalized approach \hyperref[eq: strenghened penality-based SDP]{$\textcolor{CouleurCite}{\left(\text{PmKPC}\right)^+}$} with the separating strategy \textbf{Algorithm \ref{alg: separation procedure}} results in solutions that are closer to an integer matrix.

\begin{figure}[H]
\begin{longtable}[c]{|c|cc|}
\hline
\multirow{2}{*}{$\lambda$}     & \multicolumn{2}{c|}{\textsc{Average fractionnality}}             \\ \cline{2-3} 
 &
  \multicolumn{1}{c|}{\hyperref[eq: strenghened penality-based SDP]{$\textcolor{CouleurCite}{\left(\text{PmKPC}\right)^+}$}} &
  \textbf{Algorithm \ref{alg: separation procedure}} applied to \hyperref[eq: strenghened penality-based SDP]{$\textcolor{CouleurCite}{\left(\text{PmKPC}\right)^+}$}\Tstrut \\ \hline\hline
\endhead
\hline
\endfoot
\endlastfoot
$\lambda=\exposantvanilla{0}$  & \multicolumn{1}{c|}{$2.779\exposant{-3}$}\Tstrut & $1.243\exposant{-5}$ \\
$\lambda=\exposantvanilla{-2}$ & \multicolumn{1}{c|}{$2.339\exposant{-3}$} & $1.311\exposant{-5}$ \\
$\lambda=\exposantvanilla{-4}$ & \multicolumn{1}{c|}{$1.673\exposant{-2}$} & $5.069\exposant{-5}$ \\
$\lambda=\exposantvanilla{-6}$ & \multicolumn{1}{c|}{$4.326\exposant{-2}$}\Bstrut & $8.510\exposant{-5}$ \\ \hline
\end{longtable}
\caption{Average fractionnality with the separation procedure \textbf{Algorithm \ref{alg: separation procedure}}.}\label{Fig: fractionnality with the separation procedure}
\end{figure}

\begin{wrapfigure}{r}{0.57\linewidth}
\includegraphics[width=\linewidth,keepaspectratio]{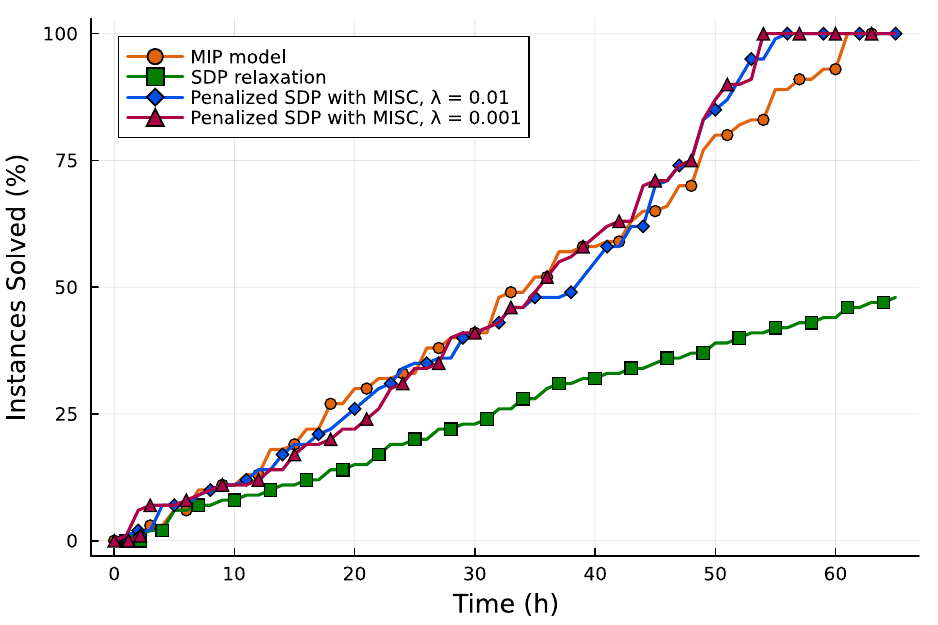}
\caption{Computation time for the studied models on the benchmark.}\label{Fig: performance profile}
\end{wrapfigure} 

\paragraph{Computing time} In terms of computing time, \cref{Fig: performance profile} plots how many instances of the benchmark are solved throughout time, depending of which model is used. The curve associated with strenghened semidefinite relaxation \hyperref[eq: strenghened semidefinite relaxation]{$\textcolor{CouleurCite}{(\text{mKPC})^{+}_{\text{SDP}}}$} shows that this model takes, on average, the longest to solve an instance, and can therefore be discarded to produce good solutions to the problem of detecting change in a time series. On the other hand, the two curves associated with the penalized models \hyperref[eq: strenghened penality-based SDP]{$\textcolor{CouleurCite}{\left(\text{PmKPC}\right)^+}$} show that the strategy of freeing ourselves from the compactness constraints allows average resolution times comparable to the MIP model. Hence the penalized approach demonstrates performance comparable to the integer model while offering the advantage of tunability, allowing for an optimal balance between compactness and a small imprecision score. Notably, for $\lambda=\exposantvanilla{-3}$, the value recommended to achieve this trade-off, the computational time per instance remains comparable to that of the MIP model.


\section{Conclusions}

In this paper, we have developed a novel semidefinite programming approach for the min-knapsack problem with compactness constraints. Our SDP model incorporates various strengthening techniques, including the derivation of valid inequalities for the semidefinite setting and a separation procedure, that together yield tighter bounds and facilitate the design of effective heuristics. In particular, the penalized approach introduces a hyperparameter $\lambda$ that allows one to balance the trade-off between solution accuracy and structural compactness, and our numerical experiments demonstrate a significant improvement in fractionality when employing a cutting plane scheme.


\section{Acknowledgments}  

This work was supported by the LabEx \textsc{Persyval}-Lab (ANR-11-LABX-0025) and by MIAI @ Grenoble Alpes (ANR-19-P3IA-0003), whose support is gratefully acknowledged.







\newpage
\appendix

\section{Useful tools of semidefinite programming}\label{app: schur complement lemma}

We recall two useful tools from the theory of semidefinite matrices that will be useful when studying reformulations of the MIP problem as a semidefinite programm:

\begin{Lemma}[Schur complement lemma]\label{lemma: Schur complement lemma}
Let $X\in\R^{n\times n}$ be the following symmetric bloc matrix
$$X=\begin{pmatrix}
   A & B^\top \\
   B & C
\end{pmatrix}$$
with $A$ invertible. Then $X\succeq 0$ if and only if $A\succeq 0$ and $C-B A^{-1}B^\top\succeq 0$.
\end{Lemma}

\begin{Theorem}[Adapted from the SDP literature, see\emph{ e.g.} theorem 3 in \citeequation{de2024integrality}]\label{thm: avec la matrice par bloc, equivalence entre X est de rang 1 et X est en 0,1}
Let $X\in\R^{n\times n}$ be a non-zero symmetric matrix such that
$$Y=\begin{pmatrix}
1 & \diag(X)^\top\\
\diag(X) & X
\end{pmatrix}\succeq 0.$$
Then $\mathrm{rank}\left(X\right)=1$ if and only if $X\in\left\lbrace 0,1\right\rbrace^{n\times n}$.
\end{Theorem}

\begin{proof}
Let us provide a direct proof of the difficult implication, different from the one of \citeequation{de2024integrality}: if $X$ has coefficients in $\left\lbrace 0,1\right\rbrace$, then $Y$ is a $\left\lbrace 0,1\right\rbrace$-positive semidefinite matrix. Such matrices have binary eigenvalues, with associated eigenvectors in $\left\lbrace 0,1\right\rbrace^{n+1}$. Ignoring zero eigenvalues, there is some $\tilde{y}_1,\dots,\tilde{y}_r\in\left\lbrace 0,1\right\rbrace^{n+1}$ such that $Y=\tilde{y}_1\tilde{y}_1^\top+\dots+\tilde{y}_r\tilde{y}_r^\top$. Because $Y_{11}=1$, we have $\left(\tilde{y}_1\tilde{y}_1^\top+\dots+\tilde{y}_r\tilde{y}_r^\top\right)_{11}=\left(\tilde{y}_1\right)_1^2+\dots+\left(\tilde{y}_r\right)_1^2=1$ which implies that only one of the $\tilde{y}_1,\dots,\tilde{y}_r$ has its first coordinate equal to one. Suppose without loss of generality that this vector is $\mathbf{y}=\tilde{y}_1$.

By definition of $Y$, for all $j\in\left\lbrace 1,\dots, n+1\right\rbrace$, we have $Y_{jj}=Y_{1j}$, thus       $$Y_{1j}=\sum^{r}_{i=1}\left(\tilde{y}_i\right)_1\left(\tilde{y}_i\right)_j=\underbrace{\mathbf{y}_1}_{=1}\mathbf{y}_j+\sum^{r}_{i=2}\underbrace{\left(\tilde{y}_i\right)_1}_{=0}\left(\tilde{y}_i\right)_j=\mathbf{y}_j=Y_{jj}$$
meaning that $\diag(Y)=\mathbf{y}$. In particular, since $\mathbf{y}$ is a $\left\lbrace 0,1\right\rbrace$-vector, we have
$$\sum^r_{i=2}\left(\tilde{y}_i\right)^2_j=Y_{jj}-\left(\tilde{y}_1\right)^2_j=\mathbf{y}_j-\underbrace{\mathbf{y}^2_j}_{=\mathbf{y}_j}=0$$
thus all $\left(\tilde{y}_i\right)_j$ are zero for $i=2,\dots,r$. Since this is true for all $j$, we deduce that all $\tilde{y}_2,\dots,\tilde{y}_r$ are zero, so $Y=\mathbf{y}\mathbf{y}^\top$, implying $\mathrm{rank}\left(Y\right)=1$. Finally, since $X$ is a submatrix of $Y$, we have $\mathrm{rank}\left(X\right)\leq\mathrm{rank}\left(Y\right)$ so $\mathrm{X}$ has rank one.
\end{proof}

\section{Non-ROAD instances}\label{subsubsec: counterexample}

Even if empirical results show that most instances have the \hyperref[conj: La conjecture]{ROAD property}, it is possible to build a family of instances of the \refequation{eq: original formulation of the mKPC} that does not verify this property. To come up with a counterexample, we start from $x^\ast$, an optimal solution of the linear relaxation \hyperref[eq: original formulation of the mKPC - LP version]{$\textcolor{CouleurCite}{\left(\text{mKPC}\right)_{\text{LP}}}$}, we set $\mathbf{X}={x^\ast}{x^\ast}^\top\,+\,\Diag({x^\ast}-{x^\ast}^2)$ and try to check if $\mathbf{X}$ verifies the constraints of the naive semidefinite relaxation. It follows easily that the knapsack constraint of \refequation{eq: naive semidefinite relaxation} is verified since
$$w^\top\diag(\mathbf{X})=\sum^n_{i=1}w_ix^\ast_i\geq q;$$
and on the other hand, since $x^\ast\in [0,1]^n$ we have $x_i^\ast-{x_i^\ast}^2\geq 0$ for all $i\in[\![n]\!]$, thus
$$\Diag\left(x^\ast-{x^\ast}^2\right)=\underbrace{x^\ast{x^\ast}^\top+ \Diag\left(x^\ast-{x^\ast}^2\right)}_{=\,\mathbf{X}} -x^\ast{x^\ast}^\top = \mathbf{X}-\diag\left(\mathbf{X}\right)\diag\left(\mathbf{X}\right)^\top\succeq 0.$$
By the \emph{Schur complement lemma} (\textbf{Lemma \ref{lemma: Schur complement lemma}} in \cref{app: schur complement lemma}), we have
$$\begin{pmatrix}
       1 & \diag\left(\mathbf{X}\right)^\top\\
       \diag\left(\mathbf{X}\right) & \mathbf{X}
   \end{pmatrix}\succeq 0$$
and the conic constraint of \refequation{eq: naive semidefinite relaxation} is also verified.

The main difficulty lies in the compacity constraints, hence we try to build instances of the \refequation{eq: original formulation of the mKPC} where the compactness constraint behaves differently depending on the formulation chosen: linear on the one hand, with a left-hand side member in $x_i+x_j-1$; and quadratic on the other hand, where the left-hand side member is $\mathbf{X}_{ij}$ ($=x_ix_j$). The key idea is to build instances whose extremal items have very large weights compared to the others, which will therefore force them to be taken via the knapsack constraint. Taking these two items alone will satisfy the knapsack constraint, and we will then have to take a minimum number of intermediate items to satisfy the compactness constraint. To create an ambiguous situation, we take $n=2m$ an even number of items with a compactness parameter $\Delta=2$: an optimal integer solution will therefore be such that exactly one in two intermediate objects is chosen. The costs of each item are chosen to be $1$, and we build the following problem:

\begin{equation}
   \left[\begin{array}{rl}
       \text{minimize} & \textbf{1}^\top x \\
       \text{subject to} & (2m+1)\left(x_1+x_{2m}\right) + \displaystyle\sum\limits_{k=2}^{2m-1}x_k\geq 4m+2\\
        & \forall i,j\in [\![n]\!],\,j-i>2\quad \left\lfloor\dfrac{j-i-1}{2}\right\rfloor\left(x_i+x_j-1\right)\leq\displaystyle\sum\limits_{k=i+1}^{j-1}x_k\\
        & x\in\left\lbrace 0,1\right\rbrace^{2m}
   \end{array}\right.\tag{$\left(\text{\textbf{CE}}\right)_m$}\label{eq: counterexample}
\end{equation}

For any $m\geq 2$, \refequation{eq: counterexample} is a integer linear program with $2m$ variables which can be seen as an instance of \refequation{eq: original formulation of the mKPC} with $2m$ items.


When we solve the linear relaxation of \refequation{eq: counterexample}, we may wonder whether the quadratic version of the compactness constraint is always verified, which would imply that \refequation{eq: counterexample} has the \hyperref[conj: La conjecture]{ROAD property}. It turns out that we can very quickly find cases where the quadratic compactness inequality is not verified. Furthermore, we can formulate the following conjecture:

\begin{Conjecture}\label{conj: etude plus approfondie du contre-exemple}
For all $m\geq 2$, there is an optimal solution $x^\circ_m$ of the linear relaxation of \emph{\refequation{eq: counterexample}} for which there exists $i,j\in\left\lbrace 1,\dots,2m\right\rbrace$ with $j-i>2$ such that
$$\left\lfloor\dfrac{j-i-1}{2}\right\rfloor \left(x^\circ_m\right)_i \left(x^\circ_m\right)_j>\sum_{k=i+1}^{j-1}\left(x_m^\circ\right)_k.$$
In particular, the matrix $\mathbf{X}:={x_m^\circ}{x_m^\circ}^\top+\Diag\left({x_m^\circ}-{x_m^\circ}^2\right)$ is not a solution of \emph{\refequation{eq: naive semidefinite relaxation}}, and therefore \emph{\refequation{eq: counterexample}} has not the \emph{\hyperref[conj: La conjecture]{ROAD property}}.
\end{Conjecture}

For example, for $m=5$, we can find $x^\circ_{5}$ an optimal solution of the linear relaxation of \refequation{eq: counterexample} which is explicitly
$$x^\circ_{5}=\begin{pmatrix}
   1 & \dfrac{3}{4} & \dfrac{119}{180} & 0 & \dfrac{17}{135} & \dfrac{251}{540} & 0 & \dfrac{107}{540} & \dfrac{11}{15} & \dfrac{11}{15}
\end{pmatrix}^\top$$
and we can note that for indexes $(i,j)=(2,9)$, the quadratic compacity constraint is not verified:
$$\left\lfloor\dfrac{9-2-1}{2}\right\rfloor \left(x_{5}^\circ\right)_2\left(x_{5}^\circ\right)_9 = 3\times\dfrac{3}{4}\times\dfrac{11}{15}=\text{\fbox{$\dfrac{33}{20}>\dfrac{29}{20}$}}=\dfrac{119}{180}+0+\dfrac{17}{135}+\dfrac{251}{540}+0+\dfrac{107}{540}=\sum^8_{k=3}\left(x_{5}^\circ\right)_k.$$

Numerically, \textbf{Conjecture \ref{conj: etude plus approfondie du contre-exemple}} has been verified for all $2\leq m\leq 800$.

Moreover, for this counter-example, we can see that the inequality in \textbf{Conjecture \ref{conj: conjecture 2 sur les bornes des problèmes}} is strict since we have, for $m=5$, that
$$\Opt\left(\text{\refequation{eq: naive semidefinite relaxation}}\right)\approx 4.42<\dfrac{14}{3}=\Opt\left(\left(\text{\textbf{CE}}\right)_5^{\text{(LP)}}\right)$$
which shows that there is a real difference in behaviour between the linear formulation of the compactness inequalities
$$\left\lfloor\dfrac{j-i-1}{\Delta}\right\rfloor\left(x_i+x_j-1\right)\leq\sum_{k=i+1}^{j-1}x_k$$
and the quadratic formulation
$$\left\lfloor\dfrac{j-i-1}{\Delta}\right\rfloor\mathbf{X}_{ij}\leq\sum_{k=i+1}^{j-1}\mathbf{X}_{kk}.$$

This leads us to conclude that it is necessary to reinforce the relaxation of the semidefinite model from \textbf{Proposition \ref{prop: formulation equivalente avec la contrainte conique}}, which is done in section \ref{sec: Strenghtening the semidefinite relaxation}. 

\section{Building the difficult instances}\label{app: Generation of instances}

Inspired from \citeequation{santini2024min}, we build our instances are constructed as follows: given an instance size $n\in\N$, two peak locations $\lambda_1$ and $\lambda_2$ are selected from the index set $ \left\lbrace 1,\dots, n\right\rbrace$. The first peak $\lambda_1$ is drawn from a normal distribution $\mathscr{N}\left(n/3, n/6\right)$, while the second peak $\lambda_2$ follows $\mathscr{N}\left(2n/3, n/6\right)$. To generate the weight distribution, we construct a histogram by sampling 5000 values for each peak, yielding two integer vectors $ w_1, w_2\in\R^n $ whose entries are sampled according to the distributions: 
$$w_1\sim\mathscr{N}(\lambda_1,n/2k)\qquad\text{and}\qquad w_2\sim \mathscr{N}(\lambda_2,n/2k)$$
where $k$ is a generation parameter chosen uniformly from $\left\lbrace 8,16,32\right\rbrace$. The final weight vector $w$ is obtained by summing the two generated vectors and normalizing them: 
$$w := \dfrac{w_1+w_2}{\left\Vert w_1+w_2\right\Vert_1}.$$

Since the case where all items have identical costs can be solved in polynomial time using dynamic programming (see Section 4.4 in \citeequation{cappello2022bayesian}), we construct the cost vector $c\in\R^n$ with no particular structure. Specifically, each coordinate $c_i$ is drawn independently from a uniform distribution over the real segment $[1,6]$. In particular, a cost cannot be below $1$, ensuring that no item is disproportionately "cheap" and thus easily selectable without significantly impacting the objective value.

The knapsack capacity $q$ is determined using a proportion parameter $p$ sampled from a uniform distribution $p\sim\mathscr{U}(0.65, 0.95)$. The capacity is then set as: 
$$q:=p \sum_{i=1}^{n} w_i.$$

Finally, the compactness parameter $\Delta$ is drawn uniformly from the set $\left\lbrace 1,2,3,4\right\rbrace$, enforcing a relatively compact solution structure.

\end{document}